\documentclass[twoside,12pt]{article}
\usepackage{amsmath}
\usepackage{amssymb}
\usepackage{mathrsfs}
\usepackage{cite}
\usepackage[usenames,dvipsnames]{color}
\definecolor{rosso}{rgb}{0.8,0,0}
\def\pier #1{{\color{rosso}#1}}
\def\takeshi #1{{\color{blue}#1}}
\def\fukao #1{{\color{red}#1}}
\let\pier\relax
\let\takeshi\relax
\let\fukao\relax

\def\interior{\mathop{\rm int}}

\title{Nonlinear diffusion equations as asymptotic limits\\
of {C}ahn--{H}illiard systems}

\author{Pierluigi Colli\\
Dipartimento di Matematica, \pier{Universit\`a di Pavia}\\
and IMATI C.N.R.\ Pavia\\
Via Ferrata~1, 27100 Pavia, Italy\\
E-mail: \texttt{pierluigi.colli@unipv.it}\\
\and \\ Takeshi Fukao\\
Department of Mathematics, Faculty of Education\\
Kyoto University of Education\\
1~Fujinomori, Fukakusa, Fushimi-ku, Kyoto~612-8522 Japan\\
E-mail: \texttt{fukao@kyokyo-u.ac.jp}}
 
\date{}

\newcommand\testopari{\sc Pierluigi Colli and Takeshi Fukao}
\newcommand\testodispari{\sc \pier{Diffusion equations as limits of {C}ahn--{H}illiard systems}}
\begin{document}
\maketitle

\begin{abstract}
An asymptotic limit \pier{of a class} of {C}ahn--{H}illiard systems \pier{is 
investigated to obtain a general} nonlinear diffusion equation. 
The target diffusion equation \pier{may reproduce a number of well-known 
model equations: {S}tefan} problem, porous media equation,
\takeshi{{H}ele-{S}haw} profile, nonlinear diffusion of singular logarithmic \pier{type},
nonlinear diffusion of {P}enrose--{F}ife type, fast diffusion equation and so on. 
Namely, by setting the suitable potential of the {C}ahn--{H}illiard systems\pier{,}
all of these problems \pier{can be obtained as limits} of 
the {C}ahn--{H}illiard \pier{related problems. Convergence results 
and error estimates are} proved.

\vspace{2mm}
\noindent \textbf{Key words:}~~{C}ahn--{H}illiard system, {S}tefan problem, 
porous media equation,
\takeshi{{H}ele-{S}haw} profile, 
fast diffusion equation. 

\vspace{2mm}
\noindent \textbf{AMS (MOS) subject clas\-si\-fi\-ca\-tion:} 35K61, 35K25, 
\pier{35B25,} 35D30, 80A22.

\end{abstract}

\section{Introduction}
\setcounter{equation}{0}

In this paper, we are interested in a discussion of the 
nonlinear diffusion problem 
\begin{gather} 
	\frac{\partial u}{\partial t} -\Delta \xi = g, 
	\quad \xi \in \beta (u) 
	\quad \mbox{in }Q:=\Omega \times (0,T), 
	\label{ND1}
	\\
	\partial _{\boldsymbol{\nu }} \xi =h
	\quad \mbox{in }\Sigma :=\Gamma \times (0,T),
	\label{ND2}
	\\
	u(0) = u_0 \quad \hbox{in}~\Omega, 
	\label{ND3}
\end{gather}
as an asymptotic limit of the following {C}ahn--{H}illiard system
\begin{gather} 
	\frac{\partial u}{\partial t} -\Delta \mu  = 0
	\quad \mbox{in }Q, 
	\label{CH1}
	\\
	\mu =-\varepsilon \Delta u + \xi +\pi _\varepsilon (u)-f, 
	\quad \xi \in \beta (u) 
	\quad \hbox{in}~Q, \label{CH2}
	\\ 
	\partial _{\boldsymbol{\nu }} \mu =\partial _{\boldsymbol{\nu }}u=0
	\quad \mbox{in }\Sigma,
	\label{CH3}
	\\
	u(0) = u_{0\varepsilon } \quad \hbox{in}~\Omega, 
	\label{CH4}
\end{gather}
as $\varepsilon \searrow 0$,
where $0<T<+\infty$ \pier{denotes} a finite time and 
$\Omega \subset \mathbb{R}^{d}$, $d=2$ or $3$, 
{\pier is} a bounded domain with smooth boundary $\Gamma $;
the symbol 
$\Delta$ stands for the {L}aplacian,
and 
$\partial _{\boldsymbol{\nu }}$ 
\pier{denotes the outward normal derivative 
on $\Gamma $.}
In the nonlinear diffusion term,  
$\beta $ is a maximal monotone 
graph and $\pi _\varepsilon $ is an anti-monotone function which 
tends to $0$ in a suitable \pier{way} as $\varepsilon \searrow 0$. 
It is \pier{well known} that the {C}ahn--{H}illiard system \eqref{CH1}--\eqref{CH4} 
is characterized by the nonlinear term $\beta +\pi_\varepsilon$, \pier{which 
represents some derivative (actually, a non-smooth $\beta $ plays as the subdifferential 
of a proper convex and lower semicontinuous function) of a multi-well function $W$.}
Usually referred as the double well potential, \pier{we can take, as a simple example,}
$W(r)=(1/4)(r^2-\varepsilon )^2$\pier{; in this 
case,} we have that $\beta (r)=r^3$ and $\pi_\varepsilon (r)=
- \varepsilon r$ for all $r \in \mathbb{R}$, \pier{$\pi_\varepsilon $ depending on} $\varepsilon >0$
\pier{(see the details of this prototype} in \cite{CH58, EZ86}). 
\pier{Actually,} the {C}ahn--{H}illiard system \eqref{CH1}--\eqref{CH4} has been
treated under various frameworks for $\beta +\pi _\varepsilon$, in particular  
\pier{for graphs $\beta $ singular with bounded domain or even non-smooth subdifferentials of bounded intervals.}
On the other hand, \pier{it is clear that 
the corresponding problem \eqref{ND1}--\eqref{ND3} may represent} various kind of nonlinear diffusion problems: {S}tefan problem, porous media equation,
\takeshi{{H}ele-{S}haw} profile, \pier{diffusion for a} singular logarithmic potential,
nonlinear diffusion of {P}enrose--{F}ife type, and 
fast diffusion equation (\pier{see the later} Examples~1--6).

The main objective of this paper is to show, \pier{for a fixed graph $\beta$ and for some known 
datum $f$ precisely related to the data $g$ and $h$ in \eqref{ND1}--\eqref{ND2}} (cf.~\eqref{case1}--\eqref{af}), the convergence \pier{of the solutions to the {C}ahn--{H}illiard system \eqref{CH1}--\eqref{CH4} 
to the respective solution of the nonlinear diffusion problem} \eqref{ND1}--\eqref{ND3}.
Namely, by \pier{performing} some asymptotic limit $\varepsilon \searrow 0$ in \eqref{CH1}--\eqref{CH4}, \pier{we naturally reaffirm the existence of
solutions} to~\eqref{ND1}--\eqref{ND3}. \pier{Of course, the initial values  
$u_{0\varepsilon } $ in \eqref{CH4} should suitably converge to the initial datum $u_0$ in \eqref{ND3}. By a similar procedure,
the {S}tefan problem with dynamic boundary conditions was obtained as asymptotic limit 
 in \cite{Fuk16}.} Moreover, \pier{as the solution to the limiting problem is 
also unique, we can prove an error estimate for the difference of solutions in suitable norms.}

A brief outline of this paper along with 
a short description of the various items is as follows\pier{.}
In Section~2, the convergence theorem is stated. 
Firstly, we set \pier{the notation that is used in the} paper. 
Next, we introduce the target problem (P) \pier{for} the nonlinear diffusion equation 
and \pier{recall the problem (P)$_\varepsilon$ for the {C}ahn--{H}illiard approximating system; we also state a} 
mathematical result for (P)$_\varepsilon$ in Proposition~2.1. 
\pier{At this point, we list and deal with various} examples for the problem~(P). 
All of these problems introduced in \pier{the} Examples~1--6 are included \pier{in the framework of Theorem~\pier{2.3}, which are stated shortly after and are focused on} the  
convergence of {C}ahn--{H}illiard systems (P)$_\varepsilon$ to the 
nonlinear diffusion problem (P). 

In Section~3, we \pier{detail the uniform estimates that will be useful to show
the convergence results}. 
In order to guarantee \pier{enough regularity for the unknowns of (P)$_\varepsilon$, 
we consider the regularized problem (P)$_{\varepsilon, \lambda}$ in which $\beta$ is replaced by its Yosida approximation $\beta_\lambda$, $\lambda >0$. 
After deriving all the estimates on (P)$_{\varepsilon, \lambda}$, 
we infer the same kind of uniform estimates and consequent regularities} for~(P)$_\varepsilon$. 

In Section~4, the proof of \pier{Theorem~\pier{2.3}} is given. The 
strategy of the proof is \pier{quite standard: by exploiting the uniform estimates, 
we pass to the limit as $\varepsilon \searrow  0$; on the other hand, let us put some
emphasis on} the monotonicity argument. 
The uniqueness for (P) is also discussed there. 

In Section~5, we prove the error estimate \pier{stated in Theorem~\pier{5.1} by applying a special} bootstrap argument. In Section~6, we \pier{can improve our results
enhancing the error estimate by Theorem~6.1}. Actually, 
under \pier{a stronger assumption for the heat source $f$, we can neglect an already required condition for the growth 
of $\beta $ and, in addition, treat a wider class of problems, in particular the ones
outlined} in Examples~5 and~6.  

A final Section~7 contains the proof of \pier{an} auxiliary proposition. A detailed index of sections and subsections is \pier{reported here.} 

\begin{itemize}
 \item[1.] Introduction
 \item[2.] \pier{Convergence} results
	\begin{itemize}
	 \item[2.1.] \pier{Notation}
	 \item[2.2.] Solution of \pier{the} {C}ahn--{H}illiard system
	 \item[2.3.] Convergence theorem
	\end{itemize}
 \item[3.] Uniform estimates
	\begin{itemize}
	 \item[3.1] Approximate problem for (P)$_\varepsilon $
	 \item[3.2] \pier{Deduction of the estimates}
	\end{itemize}
 \item[4.] Proof of the convergence theorem
 \item[5.] Error estimate
 \item[6.] Improvement of the results
 \item[7.] Appendix
\end{itemize}

\section{\pier{Convergence} results}
\setcounter{equation}{0}

In this section, we state the main results.

\subsection{\pier{Notation}}
We use the spaces $H:=L^2(\Omega )$, $V:=H^1(\Omega )$
with usual norms 
$| \cdot |_{H}$, $|\cdot |_{V}$
and inner products $(\cdot,\cdot )_{H}$, $(\cdot ,\cdot )_{V}$, respectively. 
Moreover, \pier{we introduce the space
$$W:=\{ z \in H^2(\Omega ) : \,
\partial_{\boldsymbol{\nu }}z=0 \,\hbox{ a.e.\ on~}\Gamma \}.$$
}%
The symbol $V^*$ denotes the dual space of $V$ and  
the pair $\langle \cdot ,\cdot \rangle _{V^*, V}$ \pier{stands for}
the duality pairing between $V^*$ and $V$. 
Define $m:V^* \to \mathbb{R}$ by 
\begin{equation*}
	m(z^*):=\frac{1}{|\Omega |}
	\langle z^*, 1 \rangle _{V^*,V}
	\quad \hbox{for~all~} z^* \in V^*,
	\label{mean*}
\end{equation*}
where $|\Omega |$ \pier{denotes} the volume of $\Omega $. 
If $z^* \in H$, \pier{$m(z^*)$ gives the mean value of $z^*$, i.e.,}
\begin{equation*}
	m(z^*)=\frac{1}{|\Omega |}
	\int_{\Omega }^{}z^* dx.
\label{mean}
\end{equation*}

Under these setting, we prepare a linear operator 
${\mathcal N}:D({\mathcal N}) \subseteq V^* \to V$. 
Define $D({\mathcal N}):=\{w \in V^* : m(w^*)=0 \}$\pier{;} 
then\pier{, for $w^* \in D({\mathcal N})$, we let 
$w={\mathcal N}w^*$ if} $w \in V$, $m(w)=0$ and $w$ is a solution 
of the following variational \pier{equality}
\begin{equation}
	\int_{\Omega }^{}\nabla w\cdot \nabla z dx = \langle w^*,z \rangle _{V^*,V}
	\quad \hbox{for~all~} z \in V.
	\label{n}
\end{equation}
If $w^* \in D({\mathcal N}) \cap H$, then \pier{it turns out that 
$w (={\mathcal N}w^*)$ solves the elliptic problem}
\begin{equation*} 
	\begin{cases}
	\displaystyle - \Delta w = w^* \quad \hbox{a.e.\ in~} \Omega, \vspace{1mm}\\
	\partial _{\boldsymbol{\nu }} w = 0 \quad \hbox{a.e.\ in~} \Gamma, \vspace{1mm}\\
	m(w) = 0 
	\end{cases} 
\end{equation*}
\pier{and, in particular, $w\in W$.}
We have the following property of ${\mathcal N}$: \pier{if} $v={\mathcal N}v^* $ and $w={\mathcal N}w^* $\pier{, then}
\begin{align}
	\langle w^*, {\mathcal N}v^* \rangle _{V^*,V} 
	& = \langle w^*, v \rangle _{V^*,V} = \int_{\Omega }^{}\nabla w \cdot \nabla v dx \nonumber \\
	& = \pier{\langle v^*, w \rangle _{V^*,V} 
	= \langle v^*, {\mathcal N}w^* \rangle _{V^*,V}
	\quad {\rm for~all~} v^*,w^* \in D({\mathcal N}).}
	\label{N}
\end{align}
\pier{Hence}, by defining
\begin{equation} 
	|w^*|_{V^*}^2:=\bigl |\nabla {\mathcal N} \bigl( w^*-m(w^*) \bigr )\bigr |_{H^d}^2
	+ \bigl| m(w^*) \bigr|^2 
	\quad \hbox{for~all~}w^* \in V^*,
	\label{norm}
\end{equation}
it is \pier{clear} that $|\cdot |_{V^*}$ \pier{yields} a norm of $V^*$. 

\pier{We also recall the Poincar\'e--Wirtinger inequality:
there is a positive constant $c_P$ such that 
\begin{equation} 
	 |z|_{V}^2 \le c_P
	|\nabla z|_{H^d}^2
	\quad \hbox{for~all~}z \in V \ \hbox{with}~m(z)=0.
	\label{poin}
\end{equation} 
}%

\subsection{Solution of \pier{the} {C}ahn--{H}illiard system}

In this subsection, we recall the well-known result \pier{for} the solvability for the {C}ahn--{H}illiard system. 

\pier{Let us emphasize that we term (P) the target 
problem expressed by \eqref{ND1}--\eqref{ND3}: this is an initial and  
non-homogeneous {N}eumann boundary value problem for the 
nonlinear diffusion equation \eqref{ND1}, 
where  $g : Q \to \mathbb{R}$, $h : \Sigma \to \mathbb{R}$ and 
$u_{0}: \Omega \to \mathbb{R}$, are the 
given} data. 

\pier{Moreover, for $\varepsilon >0$ we let (P)$_\varepsilon $ denote 
the {C}ahn--{H}illiard initial-boundary value problem \eqref{CH1}--\eqref{CH4}, 
in which $f : Q \to \mathbb{R}$ and  $u_{0\varepsilon }: \Omega \to \mathbb{R}$ appear as
the prescribed data and should be in some relation with $g$, $h$ and $u_0$,
as well as $\pi_\varepsilon $ has to enjoy some properties for small $\varepsilon$.} 
 
\pier{A simple remark concerning (P)$_\varepsilon $ is that, as it is  usual for 
\takeshi{Cahn--Hilliard} systems, equation~\eqref{CH1} and the second boundary 
condition in \eqref{CH3} imply consevation of the mean value for $u$, that is, $m(u(t) ) = m(u_{0\varepsilon })$ for all $t>0$. Indeed it suffices to integrate 
\eqref{CH1} by parts in space and time using \eqref{CH3} and \eqref{CH4}.
At the same time, we are interested to set the same condition of mass (or mean value) conservation for the solutions to (P).}

\pier{Thus,} we prescribe the data $g$ and $h$ such that \pier{%
\begin{equation} 
	\int_{\Omega }^{} g(t) dx + \int_{\Gamma }^{} h(t) d\Gamma =0 
	\quad \hbox{ for a.a. } t \in (0,T); 
	\label{pier4}
\end{equation}
}%
then, by simply integrating \pier{\eqref{ND1} over $\Omega $ and using \eqref{ND2},} 
we find that $(d/dt)\int_{\Omega }^{} u(t) dx =0$, whence 
\begin{equation*} 
	\frac{1}{|\Omega |} \int_{\Omega }^{} u(t) dx  = \frac{1}{|\Omega |} \int_{\Omega }^{} u_0 dx =m(u_0)=:m_0
\end{equation*}
for all $t \in [0,T]$. 
Then, we can specify $f$ acting in \pier{\eqref{CH2}} as an 
arbitrary solution of the following elliptic problem: 
\begin{equation} 
	\begin{cases}
	\displaystyle -\Delta f(t) = g(t) \quad \hbox{a.e.\ in~} \Omega, \\[2mm]
	\partial _{\boldsymbol{\nu }} f(t) = h(t) \quad \hbox{\pier{in the sense of traces on}~} \Gamma, \\
	\end{cases} 
	\label{case1}
\end{equation}
for a.a.\ $t \in (0,T)$, that is, 
\begin{align}
	\int_{\Omega }^{} \nabla f(t) \cdot \nabla z dx 
	=  \int_{\Omega }^{} g(t)  z dx + \int_{\Gamma }^{} 
	h(t)  z_\Gamma d\Gamma
	\quad \hbox{for~all}~z \in V,
	\label{af}
\end{align}
\pier{where $z_\Gamma $ denotes the trace of $z$ on $\Gamma$.} 
Throughout this paper, we assume that 
\begin{itemize}
 \item[(A1)] $\beta $ is a maximal monotone graph in 
$\mathbb{R} \times \mathbb{R}$ \pier{with effective domain $D(\beta )$ 
such that $\interior D(\beta) \ne \emptyset$,  and $\beta$} is the subdifferential
\begin{gather*}
	\beta =\partial \widehat{\beta}
\end{gather*}
of some \pier{convex and} lower semicontinuous function 
$\widehat{\beta }: \mathbb{R} \to [0,+\infty ]$ 
satisfying $\widehat{\beta }(0)=0$\pier{. This entails that
$0 \in \beta (0)$;}

 \item[(A2)] there exist \pier{two} positive constants $c_1, c_2$ such that 
\begin{equation*}
	\widehat{\beta }(r) \ge c_1 |r|^2-c_2 \quad \hbox{for~all~} r \in \mathbb{R};
\end{equation*}
 \item[(A3)] 
$\pi _\varepsilon : \mathbb{R} \to \mathbb{R}$ is a {L}ipschitz continuous function for all $\varepsilon \in (0,1]$. 
Moreover, there exist a positive constant $c_3$ and a strictly increasing continuous function 
$\sigma :[0,1] \to [0,1]$ such that $\sigma (0)=0$, $\sigma (1)=1$ and 
\begin{equation}
	\bigl| \pi_\varepsilon (0) \bigr| + |\pi '_\varepsilon |_{L^\infty (\mathbb{R})}
	\le c_3 \sigma (\varepsilon ) \quad \hbox{for~all~}\varepsilon \in (0,1];
	\label{error}
\end{equation} 
 \item[(A4)] $g \in L^2(0,T;H)$, $h \in L^2(0,T;L^2(\Gamma))$ \pier{and $g,\, h$ satisfy \eqref{pier4}. Then, we can fix a solution  $f \in L^2(0,T;V)$ of \eqref{af}};
 \item[(A5)] $u_0 \in H$ with $\widehat{\beta }(u_0) \in L^1(\Omega)$ and $m_0 \in \interior D(\beta )$. \pier{Moreover,  let  $u_{0\varepsilon } \in V$ fulfill  $m(u_{0\varepsilon })=m_0$ and 
\begin{gather} 
	|u_{0\varepsilon }|_H^2 \le \pier{c_4}, 
	\quad 
	\int_{\Omega }^{} \widehat{\beta }(u_{0\varepsilon }) dx \le \pier{c_4}, 
	\quad 
	\varepsilon |\nabla u_{0\varepsilon }|^2_{H^d} \le \pier{c_4}
	\label{apini} 
\end{gather} 
for some positive constant $c_4$ and for all $\varepsilon \in (0,1]$;  in addition, 
$u_{0\varepsilon } \to u$ strongly in $ H $  as $\varepsilon \searrow 0$.}
\end{itemize}
\pier{The existence of a family of data $\{ u_{0\varepsilon } \}$ satisfying  (A5) is checked in the Appendix.}

In order to \pier{make clear the generality of our setting}, we give here some examples in 
which assumptions (A1)--(A3) hold.

\paragraph{Example 1 \pier{[\emph{{S}tefan problem}].}}
The {S}tefan problem is \pier{a well-known model for} the mathematical description of 
the solid-liquid phase transition. \pier{A number of results is available in the literature for the Stefan problem, let us just quote e.g. \cite{Dam77, Fri68}. Using the weak formulation for this sharp interface model, as in \cite{Dam77}} one can take 
$\beta :\mathbb{R} \to \mathbb{R}$ a piecewise linear function
\pier{and $\pi _\varepsilon $ as follows:}
\begin{equation*} 
	\beta (r) =
	\begin{cases}
	k_s r & \hbox{if~} r<0, \\[2mm]
	0     & \hbox{if~} 0\le r \le L, \\[2mm]
	k_\ell (r-L) & \hbox{if}~r >L,
	\end{cases} 
	\quad 
	\pi_\varepsilon  (r) =
	\begin{cases}
	\displaystyle \varepsilon \frac{L}{2} & \hbox{if~} r<0, \\[4mm]
	\displaystyle \varepsilon \left( \frac{L}{2}-r\right) & \hbox{if~} 0\le r \le L, \\[4mm]
	\displaystyle -\frac{\varepsilon }{2}L & \hbox{if~} r >L
	\end{cases}
\end{equation*}
for all $r\in \mathbb{R}$,
where $k_s$, $k_\ell >0$ stand for the heat conductivities on the 
solid and liquid region, respectively; $L>0$ \pier{is the latent heat coefficient}. 
In this model, $u$ and $\beta (u)$ \pier{represent} the enthalpy and the 
temperature, respectively. \pier{\takeshi{One can see \cite{BP05, Fuk16, HK91}} 
and references therein about the Stefan problem and its} abstract framework. 

\paragraph{Example 2 \pier{[\emph{Porous media equation}].}} 
Let us consider the dynamics of \pier{a gas in a porous medium}. 
\pier{Let} the unknown parameter $u$ \pier{be} its density. 
The dynamics \pier{\eqref{ND1}--\eqref{ND3} is suited for this case
within the following setting of $\beta$ and corresponding} $\pi _\varepsilon \/$:
\begin{equation*}
	\beta(r)=|r|^{q-1}r, 
	\quad 
	\pi _\varepsilon (r)=-\varepsilon r
	\quad  \hbox{\pier{for} } r \in \mathbb{R}, 
\end{equation*}
\pier{with the exponent $q>1$. 
About porous media equation, there is a large amount of related work, 
\takeshi{for example, \cite{ALV84, AIMT08, AMTI06, BCS88, GV04, Hul87, Igb01, Igb02, Mar10, Yin08} and \cite{Vaz07}} so to quote a list of papers and a monograph.}

\paragraph{Example 3 \pier{[\emph{\takeshi{{H}ele-{S}haw} profile}].}}
The \takeshi{{H}ele-{S}haw} profile is characterized by 
the limiting behaviour with respect to the {S}tefan problem \pier{as the conductivities blow up.}
For simplicity, we can take $\beta $ as the inverse of the {H}eaviside graph 
${\mathcal H}$
\begin{equation}
	{\mathcal H}(r):=
	\begin{cases}
	\displaystyle 0 & \hbox{if~}\, r<0, \\[2mm]
	\displaystyle [0,1] & \hbox{if~}\, r=0, \\[2mm]
	\displaystyle 1 & \hbox{if~}\, r >0
	\end{cases} 
	\quad \hbox{\pier{for} } r \in \mathbb{R}\fukao{,} \label{pier1}
\end{equation}
that is, 
\begin{equation*}
	\beta(r)=\pier{ {\mathcal H}^{-1} (r) = {}} \partial I_{[0,1]}(r)
	\quad \hbox{\pier{for} } r \in [0,1], \quad 
	\pi _\varepsilon (r)=\varepsilon \pi (r)
	\quad \hbox{\pier{for} } r \in \mathbb{R}\fukao{,}
\end{equation*}
where $\partial I_{[0,1]}$ is 
the subdifferential of the indicator function $I_{[0,1]}$ of the interval $[0,1]$. 
The corresponding $\pi _\varepsilon $ is defined \pier{in terms of some $C^1$ function
$\pi$, strictly decreasing and vanishing} at $r=1/2$.
\pier{More details and other references can be found \takeshi{in \cite{AIMT08, AMTI06, Igb02, Kim03, PR93}}.}

\paragraph{Example 4 \pier{[\emph{Nonlinear diffusion with singular logarithmic law}].}} 
\pier{The double well potential is chosen in order that its derivative is of singular logarithmic 
type. In this case, $\beta $ is defined in an open interval, for instance $(-1,1)$, and 
it becomes singular when it approaches $-1$ and $1$.} We can take
\begin{equation*}
	\beta(r)= |r| \ln \frac{1+r}{1-r}
	\quad \hbox{\pier{for} } r \in (-1,1), \quad 
	\pi _\varepsilon (r) = - \varepsilon \, \alpha \, r
	\quad \hbox{\pier{for} } r \in \mathbb{R},
\end{equation*}
\pier{for a fixed positive constant $\alpha$. Please note that} 
the domain of $\widehat{\beta }$ is the closed interval $[-1,1]$. 
The double well structure is reproduced also in this case. \pier{Logarithmic nonlinearities 
can be found in a number of contributions devoted to the Cahn--Hilliard systems (see the recent contributions \cite{CMZ11, CMZ14, CF15a, CF15b}  and references therein.}
\medskip

\pier{The assumption (A2) plays a role for ensuring the existence of a solution to the limit 
problem, cf.~\cite{Ken90}. Actually, we can avoid it provided we replace (A4) with the stronger regularity assumption  (cf.~\eqref{pier4}) 
\begin{itemize}
 \item[(A6)] $g \in L^2(0,T;H)$, $h = 0$ \fukao{a.e.\ on $\Sigma$, $m(g(t))=0$ for a.a.\ $t\in (0,T)$. Then, let $f \in L^2(0,T;W)$ satisfy }
 \begin{align}
	\int_{\Omega }^{} \nabla f(t) \cdot \nabla z dx 
	=  \int_{\Omega }^{} g(t)  z dx 
	\quad \hbox{for~all}~z \in V. 
	\label{pier5}
\end{align}
\end{itemize}
Indeed, in view of \eqref{case1}, note that $f$ fulfils the Neumann homogeneous boundary condition and $\Delta f $ is bounded in  $L^2(0,T;L^2(\Omega))$, whence 
$f \in L^2(0,T;W)$ by elliptic regularity. Assumption~(A6) will be especially useful in Section~6 to improve the error estimate. On the other hand, the convergence result stated in Theorem~\pier{2.3} ensures the existence of a solution to the limit problem (P) as well. 
In this respect, our Theorem~\pier{2.3} turns out a generalization of \cite{Ken90}.} 

Here, we give \pier{two additional examples fitting our framework in the case when (A2) does not hold.}  

\paragraph{Example 5 \pier{[\emph{Nonlinear diffusion of {P}enrose--{F}ife type}].}} 
\pier{We take} a variation of the {S}tefan problem, written as 
\begin{equation} 
	\frac{\partial }{\partial t}\bigl (\theta +L {\mathcal H}(\theta -\theta _c) \bigr) -\Delta \left( - \frac{1}{\theta } \right) =g
	\quad \hbox{in}~Q,
	\label{slike}
\end{equation} 
where \pier{$\theta >0 $ denotes the absolute temperature, $\theta _c>0 $ is a 
critical temperature around which the phase change occurs, and the graph ${\mathcal H}$
is the same as in \eqref{pier1}}. 0
If $v$ is the selection from $\pier{\zeta (\theta) := \theta +L {\mathcal H}(\theta -\theta _c)}$, then we \pier{can} rewrite \eqref{slike} as 
\begin{equation*}
	\frac{\partial v}{\partial t} -\Delta \gamma (v) = g \quad \hbox{in}~Q\pier{,}
\end{equation*}  
\pier{where $\gamma $ is the composition of $\theta \mapsto -1/\theta$ and the inverse graph
of $\zeta.$ Since} $\gamma$ does not \pier{go} across the origin, 
we change the variable and set
\begin{equation*} 
	u = v - \theta _c, \quad \beta (u)=\gamma (u+\theta _c)+\frac{1}{\theta _c},
\end{equation*} 
\pier{in order to} match the assumption (A1). 
Note that in this case, $\pi _\varepsilon $ can be taken exactly as in \pier{the} Example~1 of the {S}tefan problem, while the assumption (A2) is not satisfied due to the 
behavior of $\beta $ as $r \to +\infty $. The \pier{limiting problem and variations of it were discussed in \cite{CL98, CLS99, CS95, CS97, FK07, KL05}.}

\paragraph{Example 6 \pier{[\emph{Fast diffusion equation}].}}
\pier{This} setting is similar to \pier{the one of} the porous media equation 
\begin{equation*}
	\beta(r)=|r|^{q-1}r, 
	\quad 
	\pi _\varepsilon (r)=-\varepsilon r,
	\quad \pier{r \in \mathbb{R},}
\end{equation*}
but with $0 < q <1$, so that there should be extinction \pier{of the solution} in a 
finite time (see, \takeshi{e.g., \cite{BC81, BCS88, Dib83, GV04, Hui07, RV02, Sab62, Vaz07}}). 
The \pier{extreme cases} for $q$ are $q=1$, which corresponds to the \pier{\emph{linear}
heat} equation, and $q=0$, which yields
\begin{equation*} 
	\beta (r)=
	\begin{cases}
	-1 & \hbox{if~} r<0, \\
	[-1,1] & \hbox{if~} r=0, \\
	1 & \hbox{if}~r >0,
	\end{cases} 
	\quad 
	\pi _\varepsilon (r)=-\varepsilon r,
	\quad \pier{r \in \mathbb{R},}
\end{equation*}
\pier{that is, a sign graph with similar behavior as ${\mathcal H}$ in \eqref{pier1}}. 
Note that \pier{whenever} $0\le q <1$, 
$\pi _\varepsilon $ turns out to act as a perturbation at \pier{infinity:} 
outside a bounded interval (whose length depends on $\varepsilon >0$) the 
potential $\widehat{\beta }+\widehat{\pi}_\varepsilon$ becomes negative. 
In all of these cases, $0$ is a local minimum and there are two absolute symmetric maxima. 
Also in this situation, assumption (A2) does not hold. 
\medskip

\pier{Now, let us go back to our general theory and formulate} an existence and uniqueness result for the problem (P)$_\varepsilon $
(see, e.g., \cite{CF15a, CV90, KN96, Kub12}).

\paragraph{Proposition~2.1.} {\it\pier{Assume \pier{either \hbox{(A1)}--\hbox{(A5)} or \hbox{(A1)}, \hbox{(A3)} with $\sigma (\varepsilon )=\varepsilon ^{1/2}$, \hbox{(A5)} and \hbox{(A6)}}. Then, for every} $\varepsilon \in (0,1]$ 
there exists a triplet $(u_\varepsilon, \mu_\varepsilon, \xi_\varepsilon )$ 
\pier{with}}
\begin{gather*}
	u_\varepsilon 
	\in H^1(0,T;V^*) \cap 
	L^\infty (0,T;V) \cap 
	L^2(0,T;W), \\
	\mu _\varepsilon 
	\in L^2(0,T;V), \quad 
	\xi _\varepsilon 
	\in L^2(0,T;H),
\end{gather*}
{\it \pier{satisfying}}
\begin{gather} 
	\bigl \langle u'_\varepsilon (t),z
	\bigr \rangle _{V^*,V}
	+ \int_{\Omega }^{} \nabla \mu _\varepsilon (t) 
	\cdot \nabla z dx
	=0 
	\quad {\it for~all~} z \in V,
	\label{weak1e} 
	\\
	\mu _{\varepsilon} (t)
	= 
	-
	\varepsilon \Delta u_{\varepsilon}(t)
	+\xi_{\varepsilon} (t)
	+\pi _\varepsilon \bigl( u_{\varepsilon} (t) \bigr) 
	-f(t) \quad {\it in}~H, 
	\label{weak2e}
\end{gather}
{\it for a.a.\ $t\in (0,T)$, \pier{and 
\begin{gather} 
	\xi _\varepsilon \in 
	\beta \bigl( u_\varepsilon \bigr)
	\quad {\it a.e.~in~} Q,
	 \label{weak3e}\\	
	 u_\varepsilon (0)=u_{0\varepsilon }
	\quad {\it a.e.~in~} \Omega.
	\label{weak4e}
\end{gather} 
}%
Moreover, there exists a positive constant $M$, independent of $\varepsilon >0$, such that}
\pier{%
\begin{gather}
	\int_{0}^{t}\bigl| u'_\varepsilon (s) \bigr|^2_{V^*} ds +
	\varepsilon \bigl| \nabla u_\varepsilon (t) \bigr|^2_{H^d} +
	\bigl| u_\varepsilon (t) \bigr|^2_{H} \le M, \label{prop1}\\
	\int_{0}^{t}\bigl| \mu _\varepsilon (s) \bigr|^2_{V} ds 
	+	\int_{0}^{t}\bigl| \xi _\varepsilon (s) \bigr|^2_{H} ds +
	\int_{0}^{t}\bigl| \varepsilon u _\varepsilon (s) \bigr|^2_{W} ds \le M, \label{prop2}
\end{gather}
}%
{\it for all $t \in [0,T]$.} 
\medskip
 
Since this type of the problem has been treated in \pier{other} papers, 
we only sketch the key \pier{points of the proof, in particular estimates \eqref{prop1}--\eqref{prop2},} in the next section.

\subsection{Convergence theorem}

In this subsection, 
we define the weak solution for the \pier{nonlinear diffusion} problem~(P). 
Then, we \pier{state the convergence} results. 

\paragraph{Definition~\pier{2.2}.} 
{\it A triplet $(u, \mu, \xi)$ \pier{with}}
\begin{gather*}
	u
	\in H^1(0,T;V^*) \cap 
	L^\infty (0,T;H), \quad
	\mu, \xi 
	\in L^2(0,T;V)
\end{gather*}
{\it is called weak solution of \hbox{(P)} if $u, \, \mu, \, \xi$ satisfy}
\begin{gather} 
	\bigl \langle u'(t),z
	\bigr \rangle _{V^*,V}
	+ \int_{\Omega }^{} \nabla \mu(t) 
	\cdot \nabla z dx
	=0 
	\quad \hbox{\it for~all~} z \in V \ \, \pier{\hbox{\it and a.a.\ $t\in (0,T)$,}	}\label{weak1} 
    \\
	\mu 
	= \xi-f, \quad 
	\xi  \in 
	\beta (u)
	\quad {\it a.e.~in~} Q, 
	\label{weak2}
	\\
	u(0)=\pier{u_0}
	\quad {\it a.e.~in~} \Omega.
	\label{weak3}
\end{gather} 

Our first result is related to the convergence of the solution of the {C}ahn--{H}illiard system {(P)$_\varepsilon $} to 
the weak solution of the nonlinear diffusion equation \hbox{(P)}. 

\paragraph{Theorem~\pier{2.3}.} 
{\it Assume \pier{either \hbox{(A1)}--\hbox{(A5)} or \hbox{(A1)}, 
\hbox{(A3)}  with $\sigma (\varepsilon ) =\varepsilon ^{1/2}$, 
\hbox{(A5)} and \hbox{(A6)}}. For each $\varepsilon \in (0,1]$, 
let $( u_\varepsilon , \mu _\varepsilon , \xi _\varepsilon )$ 
be the \pier{solution} of \hbox{(P)}$_\varepsilon $ obtained in Proposition~2.1. 
Then, there \pier{exists one triplet $( u, \mu ,\xi )$ such that}
\begin{gather*} 
	u_{\varepsilon} 
	\to u \quad \hbox{\it strongly~in~} 
	C \bigl([0,T];V^* \bigr)\ 
	\hbox{\it \pier{and} weakly~star~in~} 
	H^1(0,T;V^*) \cap L^\infty (0,T;H), 
	\\
	\varepsilon u_\varepsilon \to 0
	\quad \hbox{\it strongly~in~} 
	 L^\infty (0,T;V)\ 
	\hbox{\it \pier{and} weakly~in~} 
	L^2 (0,T;W), \\
	\pi _\varepsilon (u_\varepsilon ) \to 0
	\quad {\it strongly~in~} 
	 L^\infty (0,T;H)
\end{gather*} 
and, at least for a subsequence,
\begin{gather*}	 
	 \mu _\varepsilon \to \mu 
	\quad  {\it weakly~in~} 
	L^2 (0,T;V), \\
	\xi _{\varepsilon} 
	\to \xi \quad {\it weakly~in~} 
	L^2(0,T;H)
\end{gather*}
\pier{as $\varepsilon \searrow  0$. The triplet $( u, \mu ,\xi )$ is a weak solution 
of \hbox{(P)} and the component $u$ is uniquely determined. Moreover, 
if $\beta $ is \pier{single-valued}, then $\mu$ and $\xi $ are also unique.}}

\section{Uniform estimates}
\setcounter{equation}{0}

In this section, we obtain the uniform estimates \pier{useful to prove 
the convergence theorem. Throughout this section and the next Section~4, we 
will argue under the assumptions \hbox{(A1)}--\hbox{(A5)}; the suitable modifications 
for the other set of assumptions of Theorem~\pier{2.3} will be discussed in Section~6.}

\subsection{Approximate problem for (P)$_\varepsilon $}
In order to obtain the uniform estimates with respect to $\varepsilon \in (0,1]$,  
we need to consider \pier{a problem approximating} (P)$_\varepsilon $, \pier{and this is actually a strategy to prove Proposition~2.1}. 
Therefore, \pier{let us sketch the proof of Proposition~2.1 here. 
For each $\lambda \in (0,1]$, consider the 
problem (P)$_{\varepsilon,\lambda}$ which consists in finding the pair $(u_{\varepsilon, \lambda }, \mu _{\varepsilon, \lambda }) $ satisfying 
\begin{gather} 
	\bigl \langle u'_{\varepsilon, \lambda } (t), z
	\bigr \rangle _{V^*, V}
	+ \int_{\Omega }^{} \nabla \mu _{\varepsilon, \lambda } (t) \cdot \nabla z dx
	=0 
	\quad \hbox{for all } 
	z \in V,
	\label{weakap1} 
	\\
	\mu _{\varepsilon, \lambda }(t)
	= \lambda 
	u'_{\varepsilon, \lambda }(t)
	-
	\varepsilon \Delta u_{\varepsilon ,\lambda }(t)
	+\beta _{\lambda } \bigl( u_{\varepsilon ,\lambda } (t) \bigr)
	+\pi _\varepsilon \bigl( u_{\varepsilon ,\lambda } (t) \bigr) 
	-f(t) \quad \hbox{in}~H,
	\label{weakap2}
\end{gather}
for a.a.\ $t \in (0,T)$, and 
$$ u_{\varepsilon,\lambda}(0) =u_{0\varepsilon } \quad \hbox{in } H. $$
Here, $\beta _\lambda $ is the {Y}osida approximation of $\beta$ 
(see, e.g., \cite{Bar10, Bre73, Ken07}), that is, 
$\beta _\lambda  :\mathbb{R} \to \mathbb{R}$ is defined by}
\begin{equation*}
	\beta _\lambda  (r)
	:= \frac{1}{\lambda  } \bigl( r-J_\lambda (r) \bigr)
	:=\frac{1}{\lambda  }\bigl( r-(I+ \lambda  \beta )^{-1} (r) \bigr),
\end{equation*}
for all $ r \in \mathbb{R}$, 
where the associated $J_\lambda :\mathbb{R} \to \mathbb{R}$ is called the resolvent operator.  \pier{It is well known that $\beta _\lambda$ is the derivative of 
the {M}oreau--{Y}osida regularization $\widehat{\beta }_\lambda :\mathbb{R} \to \mathbb{R}$ 
of $\widehat{\beta }$:
\begin{equation*}
	\widehat{\beta }_{\lambda }(r)
	:=\inf_{s \in \mathbb{R}}
	\left\{ \frac{1}{2 \lambda} |r-s|^2
	+\widehat{\beta }(s) \right\} 
	= 
	\frac{1}{2 \lambda } 
	\bigl| r-J_\lambda (r) \bigr|^2 
	+ \widehat{\beta }\bigl (J_\lambda (r) \bigr )
, \quad r \in \mathbb{R}. 
\end{equation*}
The inequalities $0 \le \widehat{\beta }_\lambda (r) \le \widehat{\beta }(r)$ 
hold} for all $r \in \mathbb{R}$.

Based on \pier{available results~(cf., e.g., \cite{CF15a, KN96, Kub12}), it turns out 
that} the problem (P)$_{\varepsilon ,\lambda }$ has the solution 
$(u_{\varepsilon ,\lambda },\mu _{\varepsilon ,\lambda })$, \pier{with}  
$u_{\varepsilon ,\lambda } 
\in H^1(0,T;H) \cap C([0,T];V)
\cap L^2(0,T;W)$ and $\mu _{\varepsilon ,\lambda } \in 
L^2(0,T;V)$. 
Indeed, \eqref{weakap1}--\eqref{weakap2} 
is equivalent to an evolution equation \pier{in terms of the variable}
$v_{\varepsilon ,\lambda }:=u_{\varepsilon ,\lambda }-m_0$ in 
\pier{the} subspace $H_0:=\{ z \in H : m(z)=0\}$ of $H$. Then, 
we can \pier{adapt the results of \cite{CV90} 
for} doubly nonlinear evolution equations to show the existence.

\pier{Now, recalling the definition \pier{of ${\mathcal N}$ (see~\eqref{n}--\eqref{N})}, 
we take $z=1$ in \eqref{weakap1} and obtain $u'_{\varepsilon ,\lambda} (t) \in 
D({\mathcal N}) $; moreover, \eqref{weakap1} can be rewritten as 
\begin{equation} 
	{\mathcal N} u'_{\varepsilon ,\lambda }(t) 
	= m \bigl( \mu _{\varepsilon ,\lambda }(t) \bigr) 
	- \mu _{\varepsilon ,\lambda }(t) 
	\quad \hbox{in}~V,
	\label{Ntime}
\end{equation} 
\pier{whence \eqref{weakap2} is equivalent to} 
\begin{equation} 
	{\mathcal N} u'_{\varepsilon ,\lambda }(t) - 
	m \bigl( \mu _{\varepsilon ,\lambda }(t) \bigr) 
	+  
	\lambda 
	u'_{\varepsilon ,\lambda }(t)
	-
	\varepsilon \Delta u_{\varepsilon ,\lambda }(t)
	+\beta _{\lambda } \bigl( u_{\varepsilon ,\lambda } (t) \bigr)
	+\pi _\varepsilon \bigl( u_{\varepsilon ,\lambda } (t) \bigr) 
	= f(t) 
	\quad \hbox{in}~H,
	\label{test}
\end{equation} 
for a.a.\ $t \in (0,T)$.  As
$(d/dt)\int_{\Omega }^{} u_{\varepsilon ,\lambda }(t)dx =0$, it turns out that 
\begin{equation} 
	\pier{\frac1{|\Omega|}} \int_{\Omega }^{} u_{\varepsilon ,\lambda }(t) dx 
	= \pier{\frac1{|\Omega|}} \int_{\Omega }^{} u_{0\varepsilon }dx = m_0 
	\label{3.5a}
\end{equation} 
for a.a.\ $t \in (0,T)$, because $m(u_{0\varepsilon })=m_0$.}
	
\subsection{\pier{Deduction of the estimates}}

\pier{The proof of the convergence theorem is based on the estimates, independent 
of $\varepsilon $, for the solutions of (P)$_\varepsilon $. 
Here, we derive the useful uniform estimates on the approximating problem (P)$_{\varepsilon,\lambda}$ and, after stating and proving the series of next lemmas, we comment about the limit as $\lambda \searrow 0$.} 

\paragraph{Lemma 3.1.}
\pier{{\it There exists a positive constant $M_1$ and two values $\bar{\lambda } , \, \bar{\varepsilon  } \in (0,1], $ depending only on the data, such that}
\begin{align*}
	\int_{0}^{t}
	\bigl| u_{\varepsilon,\lambda}'(s)
	\bigr|_{V^*}^2 ds
	+ 2\lambda 
	\int_{0}^{t} 
	\bigl| u_{\varepsilon,\lambda }'(\takeshi{s})
	\bigr|_{H}^2 ds
	+ \varepsilon \bigl| \nabla u_{\varepsilon ,\lambda }(t) \bigr|_{H^d}^2 
	\\
	{} 
	+ \bigl| \widehat{\beta}_\lambda \bigl(u_{\varepsilon,\lambda }(t) \bigr) \bigr|_{L^1(\Omega )}
	+ \frac{c_1}{4} \bigl| u_{\varepsilon ,\lambda }(t) \bigr|^2_H
	\le M_1
\end{align*}
{\it for all $t\in [0,T]$, $\lambda \in (0,\bar{\lambda }]$ and  $\varepsilon \in (0,\bar{\varepsilon }]$.}}

\paragraph{Proof.} 
\pier{In order to obtain the boundedness in $L^\infty (0,T;H)$ 
for $u_{\varepsilon ,\lambda }$, here we exploit the assumption (A2). About the other set of assumptions of Theorem~2.3, we will discuss the variation of this lemma in Section 6.
We test \eqref{test} at time $s \in (0,T)$ by 
$u_{\varepsilon,\lambda}' (s) \in H$ and 
integrate the resultant over $(0,t)$. 
In view of  \eqref{N}, \eqref{norm}, and $m(u'_{\varepsilon ,\lambda }(s))=0$, we deduce
that
\begin{align*}
	\bigl( {\cal N} u_{\varepsilon ,\lambda }'(s), u_{\varepsilon ,\lambda }'(s) \fukao{\bigr)_{\! H}}
	& = \bigl \langle u_{\varepsilon ,\lambda }'(s), 
	{\cal N}u_{\varepsilon ,\lambda }'(s) \bigr \rangle _{V^*,V} \\
	& = \int_{\Omega }^{}\bigl| \nabla {\cal N} u_{\varepsilon ,\lambda }'(s) \bigr|^2dx  = \bigl|  u_{\varepsilon ,\lambda }'(s) \bigr|^2_{V^*}
\end{align*}
and $(m(\mu _{\varepsilon ,\lambda }(s)),u_{\varepsilon ,\lambda }'(s))_H=0$
for  a.a.\ $s \in (0,T)$. 
Therefore, we have}
\begin{align}
	& \int_{0}^{t} \bigl| u_{\varepsilon ,\lambda }'(s)\bigr|_{V^*}^2 ds
	+
	\lambda \int_{0}^{t}\bigl| u_{\varepsilon ,\lambda }'(s)\bigr|_{H}^2
	ds
	+ \frac{\varepsilon}{2} 
	\bigl|
	\nabla u _{\varepsilon ,\lambda }(t)
	\bigr|_{H^d}^2
	+ \int_{\Omega }^{} 
	\widehat{\beta }_\lambda 
	\bigl( u_{\varepsilon ,\lambda } (t)\bigr) dx 
	\nonumber \\
	& \quad 
	\le 
	\frac{\varepsilon}{2} 
	|
	\nabla u _{0\varepsilon }
	|_{H^d}^2
	+
	 \int_{\Omega }^{} 
	\widehat{\beta }_\lambda 
	( u_{0\varepsilon }) dx 
	+
	\int_{\Omega }^{} 
	\widehat{\pi}_\varepsilon  
	( u_{0\varepsilon} ) dx 
	- \int_{\Omega }^{} 
	\widehat{\pi}_\varepsilon  
	\bigl( u_{\varepsilon ,\lambda } (t)\bigr) dx 
	\nonumber \\
	& \quad \quad {} 
	+
	\int_{0}^{t}
	\bigl \langle u_{\varepsilon ,\lambda }'(s),f(s) \bigr \rangle _{V^*,V}
	 ds
	\label{1st} 
\end{align} 
for all $t \in [0,T]$, where 
$\widehat{\pi }_\varepsilon$ is the primitive of $\pi _\varepsilon $ defined by
$\widehat{\pi }_\varepsilon (r) :=\int_{0}^{r} \pier{\pi _\varepsilon (\rho)d\rho}$ 
for all $r \in \mathbb{R}$. 
Now, from \pier{(A2)} it follows that 
\begin{equation*} 
	\widehat{\beta }_\lambda (r)= \frac{1}{2 \lambda } 
	\bigl| r-J_\lambda (r) \bigr|^2 
	+ \widehat{\beta }\bigl (J_\lambda (r) \bigr )
	\ge  \frac{1}{{2 \bar{\lambda}}^{\vphantom K}} 
	\bigl| r-J_\lambda (r) \bigr|^2 
	+ c_1 \bigl| J_\lambda (r) \bigr|^2 -c_2
\end{equation*} 
for all $r \in \mathbb{R}$. 
Therefore, \pier{taking $\bar{\lambda }:=\min \{ 1,1/(2c_1)\}$ and} using $1/(4\bar{\lambda }) \ge c_1/2$, we have 
\begin{align}
	\int_{\Omega }^{} 
	\widehat{\beta }_\lambda 
	\bigl( u_{\varepsilon ,\lambda } (s)\bigr) dx 
	& = \frac{1}{2}\int_{\Omega }^{} 
	\widehat{\beta }_\lambda 
	\bigl( u_{\varepsilon ,\lambda } (s)\bigr) dx 
	+ \frac{1}{4 \bar{\lambda} } 
	\bigl| u_{\varepsilon ,\lambda } (s)-J_\lambda \bigl(
	u_{\varepsilon ,\lambda } (s) \bigr) \bigr|^2 _H 
	\nonumber \\
	& \quad {}
	+ \frac{c_1}{2} \bigl| J_\lambda  \bigl(
	u_{\varepsilon ,\lambda } (s) \bigr) \bigr|^2_H -\frac{c_2}{2} |\Omega |
	\nonumber \\
	& \ge 
	\frac{1}{2}\int_{\Omega }^{} 
	\widehat{\beta }_\lambda 
	\bigl( u_{\varepsilon ,\lambda } (s)\bigr) dx 
	+ \frac{c_1}{4} \bigl| 
	u_{\varepsilon ,\lambda } (s)  \bigr|^2_H -\frac{c_2}{2} |\Omega |
	\label{lem1}
\end{align}
for a.a.\ $s \in (0,T)$. 
Next, recalling the {M}aclaurin expansion and \eqref{error} of (A3), we \pier{infer that} 
\begin{equation*} 
	\bigl| \widehat{\pi }_\varepsilon (r) \bigr| \le 
	\bigl| \pi _\varepsilon (0) \bigr| |r|+\frac{1}{2} |\pi'_\varepsilon |_{L^\infty (\mathbb{R})} r^2 
	\le c_3 \, \sigma (\varepsilon ) (1+r^2) 
\end{equation*} 
for all $r \in \mathbb{R}$. \pier{Now, from (A3) we deduce that there exists $\bar{\varepsilon  } \in (0,1]$ such that $\sigma (\varepsilon ) \le c_1/(8c_3(1+|\Omega |))$
for all $\varepsilon \in (0,\bar{\varepsilon }]$.} Thus, we have 
\begin{align}
	- \int_{\Omega }^{} 
	\widehat{\pi}_\varepsilon  
	\bigl( u_{\varepsilon ,\lambda } (s)\bigr) dx 
	& 
	\pier{{}\le 
	c_3 \, \sigma (\varepsilon ) \int_{\Omega }^{} 
	\left( 1+
	\bigl| u_{\varepsilon ,\lambda } (s) \bigr| ^2
	\right)dx }
	\le 
	\frac{c_1}{8}
	\left( 1+
	\bigl| u_{\varepsilon ,\lambda } (s) \bigr| ^2_H
	\right)
	\label{lem2}
\end{align}
for a.a.\ $s \in (0,T)$. 
Moreover, using \eqref{apini} of (A5) \pier{leads to
\begin{align} 
	&\frac{\varepsilon}{2} 
	|
	\nabla u _{0\varepsilon }
	|_{H^d}^2
	+
	 \int_{\Omega }^{} 
	\widehat{\beta }_\lambda 
	( u_{0\varepsilon }) dx 
	+
	\int_{\Omega }^{} 
	\widehat{\pi}_\varepsilon  
	( u_{0\varepsilon} ) dx \nonumber \\
	&\quad \le 
	\frac{c_4}{2}
	+ 
	 \int_{\Omega }^{} 
	\widehat{\beta }
	( u_{0\varepsilon }) dx +
	c_3 \, \sigma(\varepsilon) \bigr(1+|\Omega | \bigr)
	\bigl( 1 + \bigl| u_{0\varepsilon } (s) \bigr| ^2_H \bigr)
	\le 
	\frac{3}{2} c_4
	+
	\frac{c_1}8
	\takeshi{(} 1 + \pier{c_4} \takeshi{)}.
	\label{lem3}
\end{align}}%
Thus, \pier{coll}ecting \eqref{1st}--\eqref{lem3},  \pier{with the help of 
the {Y}oung inequality we arrive at
\begin{align}
	& \frac{1}{2} \int_{0}^{t} \bigl| u_{\varepsilon ,\lambda }'(s)\bigr|_{V^*}^2 ds
	+
	\lambda \int_{0}^{t}\bigl| u_{\varepsilon ,\lambda }'(s)\bigr|_{H}^2
	ds
	+ \frac{\varepsilon}{2} 
	\bigl|
	\nabla u _{\varepsilon ,\lambda }(t)
	\bigr|_{H^d}^2
	\nonumber \\
	& \quad \quad {}
	+ \frac{1}{2} \bigl| 
	\widehat{\beta }_\lambda 
	\bigl( u_{\varepsilon ,\lambda } (t)\bigr) \bigr| _{L^1(\Omega )} 
	+ \frac{c_1}{8} 
	\bigl| u_{\varepsilon ,\lambda }(t) \bigr|_H^2 
	\nonumber \\
	& \quad 
	\le 
	\frac{c_2}{2}|\Omega | + 
	\frac{3}{2} c_4
	+
	\frac{c_1}4	
	+\frac{c_1 c_4 }{8}
	+ \frac{1}{2}
	 \takeshi{|} f \takeshi{|}_{L^2(0,T;V)}^2 
	 \label{cf}
\end{align} 
}%
for all $t \in [0,T]$. Thus, using \eqref{apini} of (A5) we see that there exists 
a positive constant $M_1$ depending only on 
\pier{$c_1$, $c_2$, $c_4$, $|\Omega |$ and $|f|_{L^2(0,T;V)}$,} independent of 
$\varepsilon \in (0,\bar{\varepsilon }]$ and $\lambda \in (0,\bar{\lambda }]$, 
such that the aforementioned estimate holds. \hfill $\Box$

\paragraph{Lemma 3.2.}
{\it There exists a positive constant $M_2$, 
independent of $\varepsilon \in (0,\bar{\varepsilon }]$ 
and $\lambda \in (0,\bar{\lambda }]$\pier{,} such that}
\begin{equation*} 
	 \int_{0}^{t} \bigl| \beta _\lambda \bigl( 
	 u_{\varepsilon ,\lambda }(s) \bigr) \bigr|_{L^1(\Omega )}^2 ds 
	 \le M_2
\end{equation*} 
{\it for all $t \in [0,T]$.}

\paragraph{Proof.} 
\pier{In view of \eqref{3.5a}, 
we have that $u_{\varepsilon ,\lambda }(s)-m_0 \in D({\cal N})$ for a.a.\ 
$s \in (0,T)$. 
Let us recall the useful inequality proved in \cite[Section~5, p.~908]{GMS09}:
thanks to (A3) and (A5) (in particular, to the facts that $m_0 $ lies in the interior of $D(\beta )$ and that $\widehat\beta(0)=0$),  it turns out that there exist two
positive constants $c_5$, $c_6$ (depending on the position of $m_0$) such that 
\begin{equation} 
	\beta_\lambda (r)(r-m_0) \ge c_5 \takeshi{\bigl|} \beta _\lambda (r) \takeshi{\bigr|} -c_6 \quad\hbox{for all } r \in \mathbb{R}.	\label{useful}
\end{equation}
Now, we can test} \eqref{test} at time $s \in (0,T)$ by $u_{\varepsilon,\lambda} 
(s) -m_0 \in D({\cal N})$.  
Then, we obtain \pier{%
\begin{align}
	& \varepsilon
	\bigl|
	\nabla u _{\varepsilon ,\lambda }(s)
	\bigr|_{H^d}^2 + 
	\bigl( \beta _\lambda \bigl(u_{\varepsilon ,\lambda }(s) \bigr ), u_{\varepsilon ,\lambda }(s)-m_0 \bigr)_{\! H}
	\nonumber \\
	& \quad \le - \bigl( {\cal N} u_{\varepsilon ,\lambda }'(s), u_{\varepsilon ,\lambda }(s)-m_0 \bigr)_{\! H}
	- \lambda \bigl( u_{\varepsilon ,\lambda }'(s) , u_{\varepsilon ,\lambda }(s)-m_0 \bigr)_{\! H}
	\nonumber \\
	& \quad 
	\quad {}- \bigl( \pi _\varepsilon  \bigl(u_{\varepsilon ,\lambda }(s) \bigr ), u_{\varepsilon ,\lambda }(s)-m_0 \bigr)_{\! H}
	+ \bigl( f(s), u_{\varepsilon ,\lambda }(s)-m_0 \bigr)_{\! H}
	\label{2nd}
\end{align}
for a.a.\ $s \in (0,T)$, because $(m(\mu _{\varepsilon ,\lambda }(s)),u_{\varepsilon ,\lambda }(s)-m_0)_H=0$.} 
Now, we know that there exists a positive constant $c_7$ such that 
$|z|_{V^*} \le c_7 |z|_H$ for all $z \in H$, therefore
\begin{align}
	- \bigl (
	{\cal N} u_{\varepsilon ,\lambda }'(s), 
	u_{\varepsilon ,\lambda }(s)-m_0
	\fukao{\bigr )_{\! H}}
	& = - \bigl \langle u_{\varepsilon ,\lambda }(s)-m_0,
	{\cal N} u_{\varepsilon ,\lambda }'(s) 
	\bigr \rangle _{V^*,V}
	\nonumber \\
	& = \bigl (
	u_{\varepsilon ,\lambda }'(s), 
	u_{\varepsilon ,\lambda }(s)-m_0
	\bigr )_{V^*} 
	\nonumber \\
	& \le c_7 \bigl| u_{\varepsilon ,\lambda }'(s) \bigr|_{V^*} 
	\Bigl( \bigl| u_{\varepsilon ,\lambda }(s) \bigr|_H+|m_0| |\Omega | \Bigr) 
	\label{lem5}
\end{align}
for a.a.\ $s \in (0,T)$. \pier{Next, we have that}
\begin{gather}
	- \lambda \bigl (
	u_{\varepsilon ,\lambda }'(s), 
	u_{\varepsilon ,\lambda }(s)-m_0
	\fukao{\bigr )_{\! H}}
	\le \lambda  \bigl| u_{\varepsilon ,\lambda }'(s) \bigr|_{H} 
	\Bigl( \bigl| u_{\varepsilon ,\lambda }(s) \bigr|_H+|m_0| |\Omega | \Bigr), 
	\label{lem6}\\
	\bigl( f(s), u_{\varepsilon ,\lambda }(s)-m_0 \bigr)_{\! H} 
	\le \bigl| f(s) \bigr|_{H} 
	\Bigl( \bigl| u_{\varepsilon ,\lambda }(s) \bigr|_H+|m_0| |\Omega | \Bigr)
	\label{lem7}
\end{gather}
for a.a.\ $s \in (0,T)$. Additionally, \pier{it is straightforward to see that} 
\begin{align}
	& - \bigl (
	\pi _\varepsilon 
	\bigl ( 
	u_{\varepsilon,\lambda } (s)
	\bigr), 
	u_{\varepsilon ,\lambda }(s)-m_0
	\fukao{\bigr )_{\! H}} \nonumber \\
	& \quad  \le \int_{\Omega }^{} \Bigl(
	\bigl| \pi _\varepsilon (u_{\varepsilon ,\lambda }(s) \bigr) - \pi _\varepsilon (0) 
	\bigr| + \bigl| \pi _\varepsilon (0) \bigr| \Bigr) 
	\bigl | 
	u_{\varepsilon ,\lambda }(s)-m_0
	\bigr| dx
	\nonumber \\
	& \quad \le \int_{\Omega }^{} 
	\Bigl( 
	| \pier{\pi _\varepsilon'} |_{L^\infty (\mathbb{R})}
	\bigl| u_{\varepsilon ,\lambda }(s) \bigr| 
	 + 
	\bigl| \pi _\varepsilon (0) \bigr| 
	\Bigr) 
	\bigl | 
	u_{\varepsilon ,\lambda }(s)-m_0
	\bigr| dx 
	\nonumber \\
	& \quad \le 
	c_3 \int_{\Omega }^{} \Bigl( \bigl| u_{\varepsilon ,\lambda }(s) \bigr| + 1 \Bigr)
	\Bigl( \bigl| u_{\varepsilon ,\lambda }(s) \bigr| + |m_0| \Bigr)dx
	\nonumber \\
	& \quad \le 
	\pier{c_8} 
	\Bigl( 1+ \bigl| u_{\varepsilon ,\lambda }(s) \bigr|_{H}^2  \Bigr)
	\label{lem8}
\end{align}
for a.a.\ $s \in (0,T)$ \pier{and for some constant $c_8 >0$ depending only on 
$c_3$, $|\Omega |$ and $|m_0|$.} 
Then, \pier{coll}ecting \eqref{2nd}--\eqref{lem8} and using \eqref{useful},
we \pier{deduce} that 
\begin{align}
	\bigl | \beta _\lambda \bigl( 
	u_{\varepsilon ,\lambda }(s) \bigr ) \bigr |_{L^1(\Omega )}
	& \le \frac{c_6}{c_5} |\Omega | 
	+ \frac{c_7}{c_5} \bigl| u_{\varepsilon ,\lambda }'(s) \bigr|_{V^*} 
	\Bigl( \bigl| u_{\varepsilon ,\lambda }(s) \bigr|_H+|m_0| |\Omega | \Bigr)
	\nonumber \\
	& \quad {} + \pier{\frac{1}{c_5} \Bigl(\lambda \bigl| u_{\varepsilon ,\lambda }'(s) 
	\bigr|_{H} +\bigl| f(s) \bigr|_{H} \Bigr) }
	\Bigl( \bigl| u_{\varepsilon ,\lambda }(s) \bigr|_H+|m_0| |\Omega | \Bigr)
	\nonumber \\
	& \quad {} + \frac{\pier{c_8}}{c_5} 
	\Bigl( 1+ \bigl| u_{\varepsilon ,\lambda }(s) \bigr|_{H}^2  \Bigr) \label{pier2}
\end{align}
for a.a.\ $s \in (0,T)$. \pier{Now, we square the sides of \eqref{pier2}, 
integrate the resultant over $(0,T)$, and take advantage of Lemma 3.1.  
Thus, we easily find a positive constant $M_2$,  
 depending only on $M_1$, $c_5$, $c_6$, $c_7$, $c_8$, $|\Omega |$, $T$, $|m_0|$ and $|f|_{L^2(0,T;H)}$, such that the assertion of the lemma follows.}
\hfill $\Box$

\paragraph{Lemma 3.3.}
{\it There exists a positive constants $M_3$,
independent of $\varepsilon \in (0,\bar{\varepsilon }]$ and $\lambda \in 
(0,\bar{\lambda }]$\pier{,} such that}

\begin{equation*}
	\int_{0}^{t} \bigl| 
	m \bigl( \mu _{\varepsilon ,\lambda } (s)\bigr) 
	\bigr|^2 ds
	\le M_3
\end{equation*}
{\it for all $t \in [0,T]$. }

\paragraph{Proof.} 
Recalling the assumption (A3), we have \pier{that}
\begin{align}
	\bigl| \pi _\varepsilon 
	\bigl( u_{\varepsilon ,\lambda } (s) \bigr) 
	\bigr|^2
	&\le \Bigl\{ \takeshi{|} \pi _\varepsilon ' \takeshi{|}_{L^\infty (\mathbb{R})} 
	\bigl|  u_{\varepsilon ,\lambda } (s) 
	\bigr| + \bigl| \pi _\varepsilon (0)
	\bigr| \Bigl\}^2 \nonumber \\
	& \le 2  c_3^2 \sigma(\varepsilon )^2 \Bigl( 1+ \bigl| u_{\varepsilon ,\lambda }(s) \bigr|^2 \Bigr) \quad \pier{\hbox{ a.e. in } \Omega,} 
	\label{piL1}
\end{align} 
for a.a.\ $s \in (0,T)$. \pier{Now, integrating \eqref{test} over $\Omega $ 
we can also exploit the Neumann homogeneous boundary condition for $u_{\varepsilon, \lambda }$ 
(hidden in the $L^2(0,T; W)$ regularity). By squaring and using \eqref{piL1},
we easily obtain} 
\begin{align*}
	\bigl| 
	m \bigl( \mu_{\varepsilon ,\lambda } (s) \bigr) 
	\bigr|^2
	& \le 
	\frac{3}{  |\Omega |^2}
	\left\{ 
	\bigl| 
	\beta _\lambda 
	\bigl( u_{\varepsilon, \lambda } (s) \bigr ) \bigr|_{L^1(\Omega )}^2
	+ |\Omega |
	\bigl| 
	\pi_\varepsilon  \bigl( u_{\varepsilon, \lambda}(s) \bigr ) 
	\bigr|_{H}^2
	+ |\Omega |
	\bigl| f(s) \bigr|^2_{H} 
	\right\} 
	\nonumber \\
	& \le \frac{3}{  |\Omega |^2}
	\left\{ 
	\bigl| 
	\beta _\lambda 
	\bigl( u_{\varepsilon, \lambda } (s) \bigr ) \bigr|_{L^1(\Omega )}^2
	+ 
	2 |\Omega | c_3^2  \Bigl( |\Omega | 
	+ \bigl| u_{\varepsilon ,\lambda }(s) \bigr|^2_H \Bigr)
	+|\Omega |
	\bigl| f(s) \bigr|^2_{H} 
	\right\} 
\end{align*} 
for a.a.\ $s \in (0,T)$, because $\sigma (\varepsilon ) \le 1$. 
Thus, by integrating \pier{over $(0,t)$,
the existence of a positive constant $M_3$ 
 depending only} on $M_1$, $M_2$, $c_3$, $|\Omega |$, 
 \pier{$T$} and $|f|_{L^2(0,T;H)}$ 
\pier{follows.}
 \hfill $\Box$ \medskip

\paragraph{Lemma 3.4.}
{\it There exists \pier{a} positive constant $M_4$,
independent of $\varepsilon \in (0,\bar{\varepsilon }]$ 
and $\lambda \in (0,\bar{\lambda }]$\pier{,} such that}
\begin{equation*}
	\int_{0}^{t} 
	\bigl| \mu _{\varepsilon,\lambda }(s) 
	\bigr|_{V}^2 ds 
	 \le M_4
\end{equation*}
{\it for all $t\in [0,T]$.}

\paragraph{Proof.} \pier{By virtue of \eqref{norm},  \eqref{poin}, \eqref{Ntime}, and 
the fact $m(u'_{\varepsilon ,\lambda }(s))=0$, we have that}
\begin{align*}
	\int_{0}^{t} \bigl| \mu _{\varepsilon,\lambda }(s) 
	\bigr|_{V}^2 ds 
	& \le 2 \int_{0}^{t} 
	\bigl| m \bigl( \mu _{\varepsilon ,\lambda }(s) \bigr) 
	\bigr| ^2_V ds
	+
	2
	\int_{0}^{t} 
	\bigl| 
	{\mathcal N} u'_{\varepsilon ,\lambda }(s) 
	\bigr| ^2_V ds 
	 \\
	 & \le 2 \int_{0}^{t} 
	\bigl| m \bigl( \mu _{\varepsilon ,\lambda }(s) \bigr) 
	\bigr| ^2_H ds
	+
	2 c_P
	\int_{0}^{t} 
	\bigl| \nabla 
	{\mathcal N} u'_{\varepsilon ,\lambda }(s) 
	\bigr| ^2_{H^d} ds 
	 \\
	 & \le 2 |\Omega | \int_{0}^{t} 
	\bigl| m \bigl( \mu _{\varepsilon ,\lambda }(s) \bigr) 
	\bigr| ^2 ds
	+
	2c_P
	\int_{0}^{t} 
	\bigl| 
	u'_{\varepsilon ,\lambda }(s) 
	\bigr| ^2_{V^*} ds 
	\le 
	M_4
\end{align*} 
for all $t \in [0,T]$, \pier{where (see Lemmas~3.1 and~3.3) 
$M_4$ is a positive constant depending only on 
$M_1$, $M_3$, $c_P$ and $|\Omega |$.} 
\hfill $\Box$

\paragraph{Lemma 3.5.}
{\it There exist two positive constants $M_5$ and $M_6$,
independent of $\varepsilon \in (0,\bar{\varepsilon }]$ and $\lambda \in (0,\bar{\lambda }]$\pier{,} such that}
\begin{gather*}
	\pier{\int_{0}^{t} 
	\bigl| \beta _\lambda \bigl( u_{\varepsilon,\lambda }(s) \bigr) 
	\bigr|_{H}^2 ds 
	 \le M_5, \quad 
	\int_{0}^{t} 
	\bigl| \varepsilon u_{\varepsilon,\lambda }(s) 
	\bigr|_{W}^2 ds 
	 \le M_6}
\end{gather*}
{\it for all $t\in [0,T]$.}

\paragraph{Proof.}
We test \eqref{weakap2} at time $s \in (0,T)$ by 
$\beta _\lambda (u_{\varepsilon,\lambda} (s))\in V$. \pier{As}
\begin{equation*}
	- \varepsilon \bigl( \Delta u_{\varepsilon ,\lambda }(s), 
	\beta _\lambda \bigl( u_{\varepsilon ,\lambda }(s)\bigr) \bigr)_{\! H} 
	= \varepsilon \int_{\Omega }^{} \beta _\lambda '\bigl( u_{\varepsilon ,\lambda }(s)\bigr)
	\bigl| \nabla u_{\varepsilon ,\lambda }(s) \bigr| ^2 dx \ge 0,
\end{equation*}
\pier{due to the monotonicity of $\beta _\lambda $, we obtain}
\begin{align*}
	& \bigl| \beta _\lambda \bigl(u_{\varepsilon ,\lambda }(s) \bigr ) \bigl|_H^2
	\nonumber \\
	& \quad \le  \bigl( \mu _{\varepsilon ,\lambda }(s)-\lambda u_{\varepsilon,\lambda  }'(s)
	-\pi _\varepsilon  \bigl(u_{\varepsilon ,\lambda }(s) \bigr )+ f(s), 
	\beta _\lambda \bigl( u_{\varepsilon ,\lambda }(s) \bigr)  \bigr)_{\! H}
	\nonumber \\
	& \quad \le 
	\frac{1}{2}\bigl| \beta _\lambda \bigl(u_{\varepsilon ,\lambda }(s) \bigr ) \bigl|_H^2
	+ 2\Bigl(  \bigl| \mu _{\varepsilon ,\lambda }(s) \bigr|_H^2
	+ \lambda^2 \bigl|  u_{\varepsilon,\lambda  }'(s) \bigr|_H^2
	+ \bigl| \pi _\varepsilon  \bigl(u_{\varepsilon ,\lambda }(s) \bigr ) \bigr|_H^2
	+ \bigl| f(s) \bigr|_H^2 \Bigr) 
\end{align*}
for a.a.\ $s \in (0,T)$. 
Now, \pier{integrating over $(0,t)$ with respect to $s$, from \eqref{piL1} and 
Lemmas~3.1 and~3.4 it follows} that 
\begin{align*}
	\int_{0}^{t} 
	\bigl| \beta _\lambda \bigl(u_{\varepsilon ,\lambda }(s) \bigr ) \bigl|_H^2
	ds 
	& \le 
	4 \int_{0}^{t} 
	\bigl| \mu _{\varepsilon ,\lambda }(s) \bigr|_H^2 
	ds 
	+ 4 \lambda^2 
	\int_{0}^{t} 
	\bigl|  u_{\varepsilon,\lambda  }'(s) \bigr|_H^2
	ds
	\nonumber \\
	& \quad {} 
	+ 8 c_3^2 \left(|\Omega|T+ \int_{0}^{t} 
	\bigl| u_{\varepsilon ,\lambda }(s) \bigr|_H^2
	ds \right) 
	+ 4 \int_{0}^{t} \bigl| f(s) \bigr|_H^2 ds 
	\le 
	M_5 
\end{align*}
for all $t \in [0,T]$, where 
$M_5$ is a positive constant depending only on 
$M_1$, $M_4$, $c_1$, $c_3$, $|\Omega |$, $T$ and $|f|_{L^2(0,T;H)}$. \pier{At last,} by comparison in \eqref{weakap2} \pier{we deduce that} 
\begin{align*}
	\int_{0}^{t} 
	\bigl| \varepsilon \Delta u_{\varepsilon ,\lambda }(s) \bigl|_{H}^2 ds
	& \le 5 
	\int_{0}^{t} 
	\bigl| \mu _{\varepsilon ,\lambda }(s)
	\bigr| ^2 ds 
	+ 5 \lambda^2 
	\int_{0}^{t} \bigl| 
	u_{\varepsilon,\lambda  }'(s) 
	\bigr| _H^2ds 
	+ 5
	\int_{0}^{t} 
	\bigl| \beta _\lambda \bigl( u_{\varepsilon ,\lambda }(s) \bigr) 
	\bigr|_H^2 ds
	\nonumber \\
	& \quad{}
	+ 10 c_3^2 
	\left(|\Omega|T+ \int_{0}^{t} 
	\bigl| u_{\varepsilon ,\lambda }(s) \bigr|_H^2
	ds \right) 
	+ 5\int_{0}^{t} 
	\bigl| f(s)
	\bigr|_H^2 ds
\end{align*}
\pier{whence, by the boundedness properties stated in Lemma~3.1 and standard elliptic regularity results, we infer
$$
	\int_{0}^{t} 
	\bigl| \varepsilon  u_{\varepsilon ,\lambda }(s) \bigl|_{W}^2 \takeshi{ds} \leq M_6
$$
for all $t \in [0,T]$, where 
$M_6$ is a positive constant having the same dependencies as $M_5$.}
\hfill $\Box$

\medskip

\pier{Here we are: by using the uniform estimates stated in 
Lemmas from~3.1 to~3.5, we can pass to the limit 
in the approximate problem (P)$_{\varepsilon ,\lambda }$ 
as $\lambda \searrow  0$, and, with the help of monotonicity arguments,  
recover a solution $(u_\varepsilon, 
\mu _\varepsilon, \xi _\varepsilon )$ to (\hbox{P})$_\varepsilon $. 
For a fixed $\varepsilon \in (0,\bar{\varepsilon }]$, the solution component 
$u_\varepsilon$ is uniquely determined and, 
if $\beta$ is single-valued, $\mu _\varepsilon$ and $ \xi _\varepsilon$ 
are unique as well (let us quote, e.g., \cite{CF15a, CV90, KN96, Kub12}
for the involved results).
Moreover, thanks to the previous lemmas, we see that, 
on the procedure of the limit as $\lambda \searrow  0$, 
all the estimates \eqref{prop1}--\eqref{prop2} hold for $(u_\varepsilon, 
\mu _\varepsilon, \xi _\varepsilon )$. These estimates turn out to be the key 
ingredient for the convergence result.}

\section{Proof of the convergence theorem}
\setcounter{equation}{0}

In this section, we prove the convergence 
theorem \pier{under the assumptions (A1)--(A5)}.

\paragraph{Proof of Theorem~\pier{2.3 (first part)}.}
Using the estimates \eqref{prop1}--\eqref{prop2}, we see that 
there exist a subsequence $\{ \varepsilon _k \}_{k \in \mathbb{N}}$, with 
$\varepsilon _k \searrow  0$ as $\pier{k\nearrow\infty} $,
and some 
limit functions $u \in H^1(0,T;V^*) \cap L^\infty (0,T;H)$,
$\mu \in L^2(0,T;V)$ and 
$\xi \in L^2(0,T;H)$ such that 
\begin{gather} 
	u_{\varepsilon_k} 
	\to u 
	\quad \hbox{weakly~star~in~} 
	H^1(0,T;V^* ) \cap L^\infty (0,T;H), 
	\label{c1}
	\\
	\varepsilon _k u_{\varepsilon_k} \to 0
	\quad \hbox{strongly~in~} 
	L^\infty (0,T;V),
	\label{c2}
	\\
	\mu _{\varepsilon_k} \to \mu 
	\quad \hbox{weakly~in~} 
	L^2(0,T;V), 
	\label{c3}
	\\
	\xi _{\varepsilon_k} 
	\to \xi  \quad \hbox{weakly~in~} 
	L^2(0,T;H)
	\label{c4}
\end{gather} 
as $\pier{k\nearrow\infty}$. 
From \eqref{c1} and the well-known {A}scoli--{A}rzela theorem
(see, e.g., \cite[Section~8, Corollary~4]{Sim87}), 
we deduce that 
\begin{equation} 
	u_{\varepsilon_k} \to u \quad \hbox{strongly~in~} 
	C\bigl( [0,T];V^* \bigr)  \pier{;}
	\label{c5}
\end{equation}
moreover, \eqref{c2} and the boundedness property in \eqref{prop2} 
imply that 
\begin{equation*} 
	\varepsilon _k u_{\varepsilon_k} \to 0
	\quad \hbox{weakly~in~} 
	L^2 (0,T;W)
	\label{c6}
\end{equation*}
as $\pier{k\nearrow\infty}$. 
With the help of the assumption (A3) (see also \eqref{piL1}), \pier{we have that 
\begin{equation*}
	\takeshi{\bigl| \pi _{\varepsilon _k}
	( u_{\varepsilon_k } ) 
	\bigr|
	\le 2\,  c_3 \, \sigma(\varepsilon_k)
	 \bigl( 1 + | u_{\varepsilon _k}| \bigr)} \quad \hbox{a.e. in } Q,
\end{equation*} 
and consequently \eqref{c1} enables us to infer that} 
\begin{equation} 
	\pi _\varepsilon (u_{\varepsilon_k} )
	\to 0 \quad \hbox{strongly~in~} 
	L^\infty (0,T;H)
	\label{c7}
\end{equation}
as $\pier{k\nearrow\infty}$. 
Now, using \eqref{c1}--\eqref{c4} and \eqref{c7}, we can pass to the limit in 
\eqref{weak1e} and \eqref{weak2e} obtaining \eqref{weak1} and the \pier{equality in} 
\eqref{weak2} for the limit functions 
$u$, $\mu $ and $\xi$. Note that the function $u$ is weakly continuous from 
$[0,T]$ to $H$, since 
\begin{equation*}
	u \in C\bigl( [0,T]; V^* \bigr) \cap L^\infty (0,T;H).
\end{equation*} 
Then, the initial condition \eqref{weak3} makes sense \pier{and follows from \eqref{weak4e} and (A5).}
Moreover, the solution component $\xi $ belongs to 
$L^2(0,T;V)$, due to $\xi=\mu+f$ \pier{and (A4)}, even thought \eqref{c4} holds true 
just in $L^2(0,T;H)$.

It remains to check that $\xi \in \beta (u)$ a.e.\ in $Q$. 
To this aim, it suffices to recall that 
\begin{equation*} 
	u_{\varepsilon _k} \to u, \quad \xi _{\varepsilon _k} 
	\to \xi 
	\quad \hbox{weakly~in~}L^2(0,T;H)
\end{equation*} 
as $\pier{k\nearrow\infty} $ and verify that 
\begin{equation} 
	\limsup_{\pier{k\nearrow\infty} } 
	\int_{0}^{T} \bigl(\xi _{\varepsilon _k}(t), u_{\varepsilon _k}(t) \bigr)_{\! H} dt 
	\le \int_{0}^{T} \bigl( \xi(t), u(t) \bigr)_{\! H} dt, 
	\label{limsup}
\end{equation} 
(cf.\ \cite[Proposition~2.2, p.~38]{Bar10}). In order to show this, 
we \pier{test \eqref{weak2} by $u_{\varepsilon_k}(t)$ and 
integrate over (0,T). Note that}
\begin{align*} 
	& \int_{0}^{T} 
	\bigl( \xi _{\varepsilon_k}(t), u_{\varepsilon _k}(t) \bigr)_{\! H} dt
	\nonumber \\
	& \quad = 
	\int_{0}^{T}\bigl( 
	\mu _{\varepsilon_k}(t) +f(t), u_{\varepsilon _k}(t) 
	\bigr)_{\! H} dt
	- \varepsilon _k \int_{0}^{T} 
	\bigl| \nabla u_{\varepsilon _k}(t) \bigr|^2 _{H^d} dt
	- \int_{0}^{T} 
	\bigl( 
	\pi _\varepsilon \bigl( 
	u_{\varepsilon_k} (t) \bigr), 
	u_{\varepsilon _k}(t) \bigr)_{\! H} dt
	\nonumber 
	\\
	& \quad \le 
	\int_{0}^{T} 
	\bigl \langle u_{\varepsilon _k}(t), \mu _{\varepsilon _k}(t)+f(t) 
	\bigr \rangle _{V^*,V} dt
	- \int_{0}^{T} \bigl( 
	\pi _\varepsilon \bigl( 
	u_{\varepsilon_k} (t)\bigr), u_{\varepsilon _k}(t)
	\fukao{\bigr )_{\! H}} dt
\end{align*}
and two terms on the last line converge. Namely, 
on account of \eqref{c3} and \eqref{c5}, we have 
\begin{align*} 
	\lim _{\pier{k\nearrow\infty} }\int_{0}^{T} 
	\bigl \langle u_{\varepsilon _k}(t), \mu _{\varepsilon _k}(t)+f(t) 
	\bigr \rangle _{V^*,V} dt
	 = \int_{0}^{T}
	\bigl \langle u(t), \mu (t)+f(t) \bigr \rangle _{V^*,V} dt \nonumber \\
	 = \int_{0}^{T} \bigl( u(t), \mu (t)+f(t)
	 \fukao{\bigr )_{\! H}} dt \pier{= \int_{0}^{T} \bigl(  \xi (t),u(t) 
	\fukao{\bigr )_{\! H}} dt,}
\end{align*}
as $\mu + f =\xi $ \fukao{a.e.\ in Q.} Moreover, \pier{\eqref{c1} and \eqref{c7}
imply} that 
\begin{align*} 
	\lim _{\pier{k\nearrow\infty}} \int_{0}^{T} \bigl( 
	\pi _\varepsilon \bigl( 
	u_{\varepsilon_k} (t)\bigr), u_{\varepsilon _k}(t)
	\fukao{\bigr )_{\! H}} dt =0.
\end{align*}
Thus, we have checked \eqref{limsup} \pier{and completely proved that the limit 
triplet $(u,\mu ,\xi )$ yields one solution of (P).}

We now verify that the solution component 
$u$ of the problem (P) is unique. 
This part \pier{follows closely} \cite[Proposition~2.1]{Ken90}. 
Assume by contradiction that there \pier{are} two solutions
\pier{$(u_i, \mu _i, \xi _i)$, $i=1,2$, and 
take the difference of the respective equations \eqref{weak1}. 
Then, we can take $z={\cal N}(u_1(s)-u_2(s))$ since 
$m(u_1(s)-u_2(s))=0$
: indeed, 
for $i=1,2$ we have
$
	\langle u_i(s),1 \rangle _{V^*,V} = 
	\langle u_0, 1 \rangle _{V^*,V} 
$
for all $s \in [0,T]$, directly from \eqref{weak1} and \eqref{weak3}. Hence, 
as $ \mu_1 - \mu_2 = \xi_1 - f - (\takeshi{\xi_2} -f)$ a.e. in $Q$ we obtain 
\begin{equation*}
	\frac{1}{2} 
	\bigl| u_1(t)-  u_2(t) 
	\bigr|_{V^*}^2 
	+ 
	\int_{0}^{t} \bigl( \xi _1(s)-\xi _2(s), u_1(s)-u_2(s) \bigr)_{\! H} ds 
	= 0
\end{equation*}
for all $t \in [0,T]$. Then, from \eqref{weak2} and the monotonicity of $\beta $ 
the second term is nonnegative, so that
\begin{equation*}
	\bigl| u_1(t)-  u_2(t) 
	\bigr|_{V^*}^2 
	= 0
\end{equation*}
for all $t \in [0,T]$ and the uniqueness property follows. 
In general, $\xi $ and $\mu $ are not uniquely determined; indeed, if $\beta$ is multivalued (as in the case of Example~3) it may happen that we 
could add to both $\xi $ and $\mu $ a function depending only on time 
and preserve the validity of \eqref{weak1} and \eqref{weak2}. 
However, if the graph $\beta $ is single-valued in its domain 
(like in Examples~1,2,4), then $\xi $ is completely determined from the inclusion in 
\eqref{weak2}, and so is $\mu $. Anyway, being $u$ unique, it turns out that the convergences 
in \eqref{c1}--\eqref{c2} and \eqref{c5}--\eqref{c7} hold for all the families 
as $\varepsilon \searrow 0$ and not only for a subesquences. Then, the assertion 
of Theorem~2.3 is completely proved.} \hfill $\Box$

\paragraph{Remark~\pier{4.1}.} The assertion of Theorem~\pier{2.3} is a consequence of our 
uniform estimates proved for all $\varepsilon \in (0,\bar{\varepsilon }]$, as 
the reader \pier{easily realizes} from the proof. The same reader can wonder whether 
the problem (P)$_\varepsilon $ has a solution for $\varepsilon \in (\bar{\varepsilon },1]$ 
\fukao{(cf.\ our statement of Proposition~2.1)}.  
The answer is yes, although the basic estimates leading to the existence of the 
solutions are different from the ones detailed here. In order to 
see which kind of \pier{approach could be followed}, we refer the reader to, e.g., 
\cite{CF15b} where a more general problem is fully investigated.

\section{Error estimate}
\setcounter{equation}{0}

We are now going to state the error estimate\pier{, still under the assumptions (A1)--(A5) but with some reinforcement. Indeed, we improve the requirement \eqref{error} in}
(A3) by prescribing that 
\begin{equation}
	\bigl| \pi_\varepsilon (0) \bigr| + |\pi '_\varepsilon |_{L^\infty (\mathbb{R})}
	\le c_3 \varepsilon ^{1/2} \quad \hbox{for~all~}\varepsilon \in (0,1],
	\label{error2}
\end{equation} 
namely $\sigma(\varepsilon ):=\varepsilon ^{1/2}$.  Moreover, {in the framework of (A5)} 
we additionally assume that 
\begin{equation}
	| u_{0\varepsilon }-u_0 |_{V^*}
	\le c_9 \varepsilon^{1/4} \quad \hbox{for~all~}\varepsilon \in (0,1],
	\label{error3}
\end{equation} 
\pier{for some positive constant $c_9$.
We observe that if we take $u_{0\varepsilon }$ exactly as in the 
choice \pier{postulated in the} Appendix, then we can find a constant $C>0$ such that
\begin{equation*}
	| u_{0\varepsilon }-u_0 |_{V^*}^2
	\le \varepsilon \langle u_0- u_{0\varepsilon }, u_{0\varepsilon } 
	\rangle_{V^*,V} \le C\varepsilon  \quad \hbox{for~all~}\varepsilon \in (0,1].
\end{equation*} 
Thus,} we have the bound \eqref{error3} even with $\varepsilon ^{1/2}$.

\paragraph{Theorem~\pier{5.1}.} 
{\it Assume \hbox{(A1)}--\hbox{(A5)} with \eqref{error2} and \eqref{error3}. 
For $\varepsilon \in (0,\bar{\varepsilon }]$, let 
$(u_\varepsilon ,\mu_\varepsilon ,\xi _\varepsilon )$ be a solution of 
problem \hbox{(P)}$_\varepsilon $ and let $(u,\mu ,\xi )$ \pier{solve} problem \hbox{(P)}. Then, there exists a constant $C^*>0$, depending only on the data, such that}
\begin{equation} 
	| u_\varepsilon -  u 
	|_{C([0,T];V^*)}^2 
	+ 
	\int_{0}^{T} \bigl( \xi _\varepsilon (s)-\xi (s), u_\varepsilon (s)-u(s) \bigr)_{\! H} ds 
	\le C^* \varepsilon ^{1/3}
	\label{error0}
\end{equation} 
{\it for all $\varepsilon \in (0,\bar{\varepsilon }]$.}

\bigskip

Observe that the elements $u_\varepsilon $ and $u$ appearing in the error estimate 
\eqref{error0} are uniquely determined.

\paragraph{Proof \pier{of Theorem~5.1}.} 
\pier{For $\varepsilon \in (0,\bar{\varepsilon }]$, let us use  
\eqref{weakap2} at the level of the $\lambda$-approx\-imation. We rewrite  
\eqref{weakap2} as 
\begin{equation*} 
	-\varepsilon \Delta u_{\varepsilon, \lambda} 
	+ \beta_\lambda (u_{\varepsilon, \lambda}) = 
	-\lambda u_{\varepsilon ,\lambda }'
	+ \mu _{\varepsilon, \lambda} 
	+ f 
	-\pi_\varepsilon (u_{\varepsilon, \lambda } )
	\quad \hbox{a.e.\ in~}Q,
\end{equation*} 
and we test it by $-\varepsilon ^\alpha \Delta u_{\varepsilon ,\lambda }$ with $\alpha >0$. 
Then, since  
\begin{align*}
	\int_{0}^{T} 
	\bigl( \beta_\lambda \bigl( u_{\varepsilon ,\lambda }(s) \bigr) , 
	-\varepsilon ^\alpha \Delta u_{\varepsilon ,\lambda }(s) \bigr)_{\! H} ds
	& =\varepsilon ^\alpha 
	\int_{0}^{T} \! \! 
	\int_{\Omega }^{}
	\beta _\lambda ' \bigl( u_{\varepsilon ,\lambda } (s) \bigr) 
	\bigl|
	\nabla u_{\varepsilon ,\lambda }(s) 
	\bigr|^2 
	dx
	ds
	\ge 0.
\end{align*}
and 
\begin{align*}
	\int_{0}^{T} 
	\bigl( -\lambda u_{\varepsilon ,\lambda }'(s), 
	-\varepsilon ^\alpha \Delta u_{\varepsilon ,\lambda }(s) \bigr)_{\! H} ds
	& = -\lambda \varepsilon ^\alpha 
	\int_{0}^{T}
	\frac{1}{2}\frac{d}{ds}
	\int_{\Omega }^{} \bigl| \nabla u_{\varepsilon ,\lambda }(s) \bigl|^2 dx 
	ds \\
	& = -\lambda \varepsilon ^\alpha 
	\left(  \frac{1}{2} 
	\int_{\Omega }^{} \bigl| \nabla u_{\varepsilon ,\lambda }(T) \bigl|^2 dx 
	-
	 \frac{1}{2} 
	\int_{\Omega }^{} | \nabla u_{0\varepsilon} |^2 dx
	\right) \\
	& \le \frac{\lambda \varepsilon ^\alpha }{2} 
	\int_{\Omega }^{} | \nabla u_{0\varepsilon} |^2 dx = 
	\frac{\lambda \varepsilon ^\alpha }{2} 
	| \nabla u_{0\varepsilon} |^2_{H^d},
\end{align*}
using the Young inequality, we obtain 
\begin{align*}
	& \varepsilon ^{1+\alpha }\int_{0}^{T} 
	\bigl| \Delta u_{\varepsilon ,\lambda }(s) \bigr|^2_H  ds \\
	& \quad \le 
	\int_{0}^{T} 
	\bigl( \mu _{\varepsilon ,\lambda } (s) + f(s) 
	-\pi _\varepsilon \bigl( u_{\varepsilon ,\lambda } (s) \bigr),  
	-\varepsilon ^\alpha \Delta u_{\varepsilon ,\lambda }(s) \bigr)_{\! H}
	ds + \frac{\lambda \varepsilon ^\alpha}{2} 
	| \nabla u_{0\varepsilon} |^2 _{H^d}
	\\
	& \quad =  
	\int_{0}^{T} \bigl( \nabla  \bigl( \mu _{\varepsilon ,\lambda } (s)+f(s) \bigr) , 
	 \varepsilon ^\alpha \nabla u_{\varepsilon ,\lambda }(s) \bigr)_{\! H}
	ds
	\nonumber \\
	& \quad \quad {}- \int_{0}^{T} \bigl( \nabla \pi _\varepsilon 
	\bigl( u_{\varepsilon ,\lambda } (s) \bigr), 
	\varepsilon ^\alpha  \nabla u_{\varepsilon ,\lambda }(s) \bigr)_{\! H}
	ds
	+
	\frac{\lambda \varepsilon ^\alpha}{2} 
	| \nabla u_{0\varepsilon} |^2 _{H^d} \\
	& \quad \le \frac{1}{2} 
	\int_{0}^{T} \bigl| \nabla \bigl( 
	\mu _{\varepsilon ,\lambda } (s) +f(s) \bigr) \bigl| _{H^d}^2 ds 
	+ \frac{1}{2} \int_{0}^{T} 
	\bigl| \pi _\varepsilon' \bigl( u_{\varepsilon ,\lambda } (s) \bigr) 
	\nabla u_{\varepsilon ,\lambda } (s)
	\bigr| _{H^d}^2 ds 
	\nonumber \\
	& \quad \quad {}
	+ \varepsilon ^{2\alpha} \int_{0}^{T} 
	\bigl| \nabla u_{\varepsilon ,\lambda }(s) \bigr|_{H^d}^2 
	ds +
	\frac{\lambda \varepsilon ^\alpha}{2} 
	| \nabla u_{0\varepsilon} |^2_{H^d}.
\end{align*}
Then, in view of \eqref{error2} we find out that 
\begin{align}
	\varepsilon ^{1+\alpha } \int_{0}^{T} 
	\bigl| \Delta u_{\varepsilon ,\lambda }(s) \bigr|^2_H  ds 
	& \le
	\takeshi{|} \mu _{\varepsilon ,\lambda } \takeshi{|}_{L^2 (0,T;V)}^2
	+
	\takeshi{|} f \takeshi{|}_{L^2 (0,T;V)}^2
	+ \frac{1}{2}  c_3^2 T \varepsilon
	\takeshi{|} 
	\nabla u_{\varepsilon ,\lambda }
	\takeshi{|} _{L^\infty (0,T;H^d)}^2  
	\nonumber \\
	& {} \quad 
	+ \varepsilon ^{2\alpha} \int_{0}^{T} 
	\bigl| \nabla u_{\varepsilon ,\lambda }(s) \bigr|_{H^d}^2 
	ds +
	\frac{\lambda \varepsilon ^\alpha}{2} 
	| \nabla u_{0\varepsilon} |^2_{H^d}.
	\label{boot}
\end{align}
Now,} we know that $u_{\varepsilon ,\lambda } \in L^2(0,T;W)$, therefore
\begin{align}
	\varepsilon ^{2\alpha} \int_{0}^{T} 
	\bigl| \nabla u_{\varepsilon ,\lambda }(s) \bigr|_{H^d}^2 
	ds
	& = \varepsilon ^{2\alpha }
	\int_{0}^{T} \bigl(  u_{\varepsilon ,\lambda }(s), -\Delta u_{\varepsilon ,\lambda }(s)
	\bigr) _{\takeshi{H}} ds 
	\nonumber \\
	& \le \frac{1}{2} \int_{0}^{T} \bigl| u_{\varepsilon ,\lambda }(s) 
	\bigr|_{H}^2 ds + \frac{\varepsilon ^{4\alpha }}{2} 
	\int_{0}^{T} \bigl| \Delta u_{\varepsilon ,\lambda }(s) 
	\bigr|_{H}^2 ds.
	\label{boot2}
\end{align}
If, we devise $\varepsilon ^{4\alpha} \le \varepsilon ^{1+\alpha} $, that is 
$\alpha = 1/3$, \pier{then we recall Lemmas 3.1 and 3.4 and point out that \eqref{boot} and \eqref{boot2} imply}
\begin{align}
	\varepsilon ^{4/3} \int_{0}^{T} 
	\bigl| \Delta u_{\varepsilon ,\lambda }(s) \bigr|^2_H ds 
	& \le
	\left( 2M_4 
	+
	2 \takeshi{|} f \takeshi{|} _{L^2(0,T;V)}^2 
	+ c_3^2 T M_1
	+ \frac{4TM_1}{c_1}
	\right) 
	+
	\lambda \varepsilon ^\alpha
	| \nabla u_{0\varepsilon} |^2_{H^d}
	\nonumber \\
	& =: \pier{C^*_1}+
	\lambda \varepsilon ^\alpha
	| \nabla u_{0\varepsilon} |^2_{H^d}.
	\label{boot3}
\end{align}
\pier{The estimate \eqref{boot3} works for the approximate solution $u_{\varepsilon ,\lambda }$. However, thanks to Proposition~2.1 and the proof of the convergence theorem, which 
was treated in the previous section, 
$(u_{\varepsilon ,\lambda },\mu_{\varepsilon, \lambda}, \beta _\lambda (u_{\varepsilon, \lambda}))$ converges to 
$(u_\varepsilon,\mu _\varepsilon ,\xi _\varepsilon )$  as $\lambda \searrow 0$ and, passing to the limit in \eqref{boot3}, we recover the key estimate} 
\begin{equation}
	\varepsilon ^{4/3} \int_{0}^{T} 
	\bigl| \Delta u_{\varepsilon }(s) \bigr|^2_H  ds 
	\le \pier{C^*_1}.
	\label{boot4}
\end{equation}
Next, take the difference between \eqref{weak1e} and \eqref{weak1}, then \pier{we have}
\begin{equation*} 
	\bigl \langle u'_\varepsilon (s)-u'(s),z
	\bigr \rangle _{V^*,V}
	+ \int_{\Omega }^{} \nabla \bigl( \mu_\varepsilon (s)-\mu (s) \bigr)  
	\cdot \nabla z dx
	=0 
	\quad \hbox{for~all~} z \in V,
\end{equation*}
for a.a.\ $s \in (0,T)$\pier{. We can choose $z={\cal N}(u_\varepsilon (s)-u(s))$ 
because (cf.~(A5)) $m(u_\varepsilon (s)-u(s))=m(u_{0\varepsilon}) -m (u_0 )=0$;} 
then we obtain
\begin{equation*} 
	\frac{1}{2}\frac{d}{ds}\bigl| u_\varepsilon (s)-u(s)
	\bigr|_{V^*}^2
	+ \bigl \langle u_\varepsilon (s)-u(s), \mu_\varepsilon (s)-\mu (s) 
	\bigr \rangle _{V^*,V} 
	=0 
\end{equation*}
for a.a.\ $s \in (0,T)$. 
Now, 
\pier{by integrating over $(0,t)$,
owing to  \eqref{weak2e} and \eqref{weak2} we infer 
the following equality:}
\begin{align} 
	& \frac{1}{2} \bigl| u_\varepsilon (t)-u(t) \bigr|_{V^*}^2 
	+ \int_{0}^{T} 
	\bigl( u_\varepsilon (s)-u(s), \xi _\varepsilon (s)-\xi (s) \bigr)_{\! H} ds 
	\nonumber \\
	& \quad = \frac{1}{2} |u_{0\varepsilon }-u_0|_{V^*}^2 + 
	\varepsilon \int_{0}^{T}
	\bigl( \Delta u_\varepsilon (s),u_\varepsilon (s)-u(s) \bigr)_{\! H} ds
	\nonumber \\
	& \quad \quad {}- \int_{0}^{T} \bigl( \pi _\varepsilon \bigl( u_\varepsilon (s) \bigr), 
	u_\varepsilon (s)-u(s) \bigr)_{\! H} ds ,
	\label{key}
\end{align}
\pier{holding for all $t\in [0,T].$ 
Note that, by using Lemma 3.1 and \eqref{boot4}, one can find} a positive constant $C_2^*$, depending only on  $c_1$, $M_1$, $C_1^*$ and $T$, such that
\begin{align} 
	& 
	\varepsilon \int_{0}^{T}
	\bigl( \Delta u_\varepsilon (s),u_\varepsilon (s)-u(s) \bigr)_{\! H} ds
	\nonumber \\
	& \quad \le 
	\varepsilon ^{1/3} 
	\left\{ \varepsilon ^{4/3} \int_{0}^{T} \bigl| \Delta u_\varepsilon (s) \bigr|_H^2 ds \right\}^{1/2} 
	\left\{ \int_{0}^{T} \bigl| u_\varepsilon (s)-u(s) \bigr|_H^2 ds \right\}^{1/2} 
	 \le C_2^*  \varepsilon ^{1/3}.
	\label{key1}
\end{align}
Next, \pier{on account of} Lemma 3.1 and \eqref{error2}, there is a 
positive constant $C_3^*$, depending only on 
$c_1$, $c_3$, $M_1$, $|\Omega |$ and $T$, such that
\begin{align} 
	& 
	- \int_{0}^{T} \bigl( \pi _\varepsilon \bigl( u_\varepsilon (s) \bigr), 
	u_\varepsilon (s)-u(s) \bigr)_{\! H} ds 
	\nonumber \\
	& \quad \le 
	|\pi _\varepsilon '|_{L^\infty (\mathbb{R})} 
	\left\{ \int_{0}^{T} \bigl| u_\varepsilon (s) \bigr|_H^2 ds 
	\right\}^{1/2} 
	\left\{ \int_{0}^{T} \bigl| u_\varepsilon (s)-u(s) \bigr|_H^2 ds \right\}^{1/2} 
	\nonumber \\
	& \quad \quad {}
	+ \bigl| \pi _\varepsilon (0) \bigr| \, 
	\pier{\bigl(|\Omega | T \bigr)^{1/2}} 
	\left\{ \int_{0}^{T} \bigl| u_\varepsilon (s)-u(s) \bigr|_H^2 ds \right\}^{1/2} 
	\le C_3^* \varepsilon ^{1/2}.  
	\label{key2}
\end{align}
Thus, \pier{collecting} \eqref{key}--\eqref{key2} and \pier{applying 
\eqref{error3} lead to} 
\begin{align*} 
	| u_\varepsilon -u |_{C([0,T];V^*)}^2
	& \le \pier{c_9}^2 \varepsilon^{1/2}+2C_2^* \varepsilon ^{1/3} + 
	 2C_3^* \varepsilon ^{1/2},
\end{align*}
and 
\begin{equation*} 
	\int_{0}^{T} \bigl( u_\varepsilon (s)-u(s), \xi _\varepsilon (s)-\xi(s) \bigr)_{\! H} ds 
	\le \frac{1}{2} 
	\pier{c_9}^2 \varepsilon ^{1/2}+C_2^* \varepsilon ^{1/3} + 
	C_3^* \varepsilon ^{1/2}, 
\end{equation*}
that is, there exists $C^*>0$ such that the error estimate \eqref{error0} holds. 
\hfill $\Box$

\paragraph{Remark~\pier{5.2}.} If $\beta $ is {L}ipschitz continuous, \pier{then  
\eqref{error0} gives us the additional information
\begin{align*}
\int_{0}^{T} 
	\bigl| \xi _\varepsilon (s)-\xi (s) \bigr|^2_H ds 
	& \le C_\beta \int_{0}^{T} 
	\bigl( \xi _\varepsilon (s)-\xi (s), u_\varepsilon (s)-u(s) \bigr)_{H} ds 
	\le C_\beta C^* \varepsilon ^{1/3},
\end{align*}
due to the monotonicity of $\beta $. Here, $C_\beta $ denotes a {L}ipschitz constant for} $\beta $. 

\section{Improvement of the results}
\setcounter{equation}{0}
\pier{Throughout this section, we assume that \hbox{(A1)}, \hbox{(A3)} with $\sigma (\varepsilon )=\varepsilon ^{1/2}$ (that is, \eqref{error2}), \hbox{(A5)} and \hbox{(A6)} hold. Let us note that the assumption (A6) implies that} 
\begin{equation*} 
	\begin{cases}
	\displaystyle -\Delta f(t) = g(t) \quad \hbox{a.e.\ in~} \Omega, \\[2mm]
	\partial _{\boldsymbol{\nu }} f(t) = 0 \quad \hbox{a.e.\ in~} \Gamma, \\
	\end{cases} 
\end{equation*}
\pier{for a.a. $t\in (0,T)$ (cf. \eqref{case1}). We point out that (A2) is no 
longer in use and (A4) is covered by (A6).}

\paragraph{Proof of Theorem~\pier{2.3 (final part).}} 
\pier{Referring to the first part of the proof given in Section~3, here we 
have to modify  the proof of Lemma~3.1. We} take $z:=u_{\varepsilon ,\lambda }(s)$ in \eqref{weakap1} \pier{obtaining} 
\begin{equation} 
	\bigl \langle u'_{\varepsilon, \lambda } (s), u_{\varepsilon, \lambda } (s)
	\bigr \rangle _{V^*, V}
	+ \int_{\Omega }^{} \nabla \mu _{\varepsilon,\lambda } (s)
	\cdot \nabla u_{\varepsilon, \lambda } (s) dx
	=0 
	\label{6-1}
\end{equation}
for a.a.\  $s \in (0,T)$\pier{; on the other hand, 
we test \eqref{weakap2} by $-\Delta u_{\varepsilon ,\lambda }(s)$ and 
exploit~(A6) to deduce that} 
\begin{align}
	& \int_{\Omega }^{} 
	\nabla  \mu _{\varepsilon,\lambda } (s)\cdot \nabla  u_{\varepsilon ,\lambda }(s) 
	dx
	\nonumber \\
	& \quad = \frac{\lambda}{2} \frac{d}{ds}
	\bigl| \nabla 
	u_{\varepsilon ,\lambda }(s)
	\bigr|_{H^d}^2
	+\varepsilon \bigl| \Delta u_{\varepsilon ,\lambda }(s) \bigr| ^2_H
	+\int_{\Omega }^{} 
	\beta _{\lambda }' \bigl( u_{\varepsilon ,\lambda } (s) \bigr)
	\bigl| \nabla u_{\varepsilon, \lambda } (s) \bigr|^2 dx 
	\nonumber \\
	& \quad \quad {} 
	- \bigr( \pi _\varepsilon \bigl( u_{\varepsilon ,\lambda } (\pier{s}) \bigr), \Delta u_{\varepsilon ,\lambda }(s) \bigr)_{\! H}  
	+ \bigl( \Delta f(\pier{s}),u_{\varepsilon ,\lambda }(s)\bigr)_{\! H} 
	\label{6-2}
\end{align}
for a.a.\ $s \in (0,T)$.
\pier{By virtue of (A3) and the Young inequality, we have} that 
\begin{align}
	 \bigr( \pi _\varepsilon \bigl( u_{\varepsilon ,\lambda } (\pier{s}) \bigr), \Delta u_{\varepsilon ,\lambda }(s) \bigr)_{\! H}  
	& \le \int_{\Omega }^{} 
	\Bigl( |\pi _\varepsilon' |_{L^\infty (\mathbb{R})} 
	\bigl| u_{\varepsilon ,\lambda }(s) \bigr| + 
	\bigl| \pi _\varepsilon (0) \bigr| \Bigr) 
	\bigl| \Delta u_{\varepsilon ,\lambda } (s) \bigr| dx 
	\nonumber \\
	& \le \int_{\Omega }^{} c_3 \varepsilon ^{1/2}
	\Bigl( 
	\bigl| u_{\varepsilon ,\lambda }(s) \bigr| + 1 \Bigr) 
	\bigl| \Delta u_{\varepsilon ,\lambda } (s) \bigr| dx 
	\nonumber \\
	& \le \frac{\varepsilon }{2} 
	\int_{\Omega }^{}\bigl| \Delta u_{\varepsilon ,\lambda } (s) \bigr|^2 dx 
	 + \frac{c_3^2}{2}\int_{\Omega }^{}
	 \Bigl( 
	\bigl| u_{\varepsilon ,\lambda }(s) \bigr| + 1 \Bigr) ^2 dx 
	\nonumber \\
	& \le \frac{\varepsilon }{2} 
	 \bigl| \Delta u_{\varepsilon ,\lambda }(s) \bigr| ^2_H
	 + c_3^2 \Bigl( 
	\bigl| u_{\varepsilon ,\lambda }(s) \bigr|_H^2 + |\Omega | \Bigr) 
	\label{6-3}
\end{align}
and 
\begin{equation}
	\bigl( \Delta f(s),u_{\varepsilon ,\lambda }(s)\bigr)_{\! H} 
	\le \frac{1}{2} \bigl| \Delta f(s) \bigr|_{H}^2 + 
	\frac{1}{2} \bigl| u_{\varepsilon ,\lambda }(s) \bigr|_{H}^2 
	\label{6-4}
\end{equation}
for a.a.\ $s \in (0,T)$. 
Therefore, \pier{by combining \eqref{6-1}--\eqref{6-4} and integrating 
over $(0,t)$} with respect to $s$, we \pier{infer} that
\begin{align}
	& \frac{1}{2}
	\bigl| 
	u_{\varepsilon ,\lambda }(t)
	\bigr|_{H}^2
	+ \frac{\lambda}{2}
	\bigl| \nabla 
	u_{\varepsilon ,\lambda }(t)
	\bigr|_{H^d}^2
	+ \frac{\varepsilon }{2} 
	\int_{0}^{t}
	\bigl| \Delta u_{\varepsilon ,\lambda }(s) 
	\bigr|^2_H ds
	\nonumber \\
	& \quad \le \frac{1}{2} 
	| u_{0\varepsilon }
	|_{H}^2
	+ \frac{\pier{\lambda}}{2} 
	| \nabla u_{0\varepsilon }
	|_{H^d}^2
	+ c_3^2 \int_{0}^{t}\Bigl( 
	\bigl| u_{\varepsilon ,\lambda }(s) \bigr|_H^2 + |\Omega | \Bigr) ds
	\nonumber \\
	& \quad  \quad {} 
	+ \frac{1}{2} \int_{0}^{t} \bigl| \Delta f(s) \bigr|_{H}^2 ds+ 
	\frac{1}{2} \int_{0}^{t} \bigl| u_{\varepsilon ,\lambda }(s) \bigr|_{H}^2 ds
	\label{6-5}
\end{align}
for all $t \in [0,T]$, 
\pier{thanks to the monotonicity of $\beta _\lambda $. At this point, we let $\lambda\leq \varepsilon $ (which is always possible since $\lambda$ is going to $0$ before $\varepsilon$)
and recall \eqref{apini} that entails 
$$ \frac{1}{2} 
	| u_{0\varepsilon }
	|_{H}^2
	+ \frac{\pier{\lambda}}{2} 
	| \nabla u_{0\varepsilon }
	|_{H^d}^2
    \le c_4 .
$$
Then, we  apply the {G}ronwall inequality and obtain 
\begin{align}
	& \bigl| 
	u_{\varepsilon ,\lambda }(t)
	\bigr|_{H}^2 + \lambda
	\bigl| \nabla 
	u_{\varepsilon ,\lambda }(t)
	\bigr|_{H^d}^2
	+ \varepsilon  
	\int_{0}^{t}
	\bigl| \Delta u_{\varepsilon ,\lambda }(s) 
	\bigr|^2_H ds
	\nonumber
	\\
	& \quad \le 
	\left\{ c_4
	+ |\Delta f|_{L^2(0,T;H)}^2+ 2c_3^2|\Omega |T
	\right\} \exp \fukao{\bigl\{}(2c_3^2+1)T \fukao{ \bigr\} }
	\label{6-6}
\end{align}
for all $t \in [0,T]$. Hence, going back to the proof of Lemma 3.1, let us point out 
the inequalities \eqref{1st} and \eqref{lem1}: here we do not use \eqref{lem1} and  
from \eqref{1st} we arrive at 
\begin{align}
	& \frac{1}{2} \int_{0}^{t} \bigl| u_{\varepsilon ,\lambda }'(s)\bigr|_{V^*}^2 ds
	+
	\lambda \int_{0}^{t}\bigl| u_{\varepsilon ,\lambda }'(s)\bigr|_{H}^2
	ds
	\nonumber \\
	& \quad {}+ \frac{\varepsilon}{2} 
	\bigl|
	\nabla u _{\varepsilon ,\lambda }(t)
	\bigr|_{H^d}^2
	+ \bigl| 
	\widehat{\beta }_\lambda 
	\bigl( u_{\varepsilon ,\lambda } (t)\bigr) \bigr| _{L^1(\Omega )} 
	\le 
	\frac{3}{2} c_4
	+
	\frac{c_1}4	
	+\frac{c_1 c_4 }{8}
	+ \frac{1}{2}
	 \bigl| f\bigr|_{L^2(0,T;V)}^2 
	 \label{pier3}
\end{align} 
for all $t \in [0,T]$, which replaces \eqref{cf}. Now,  by adding
\eqref{6-6} and \eqref{pier3} we obtain the useful bound to
continue with other lemmas and end the proof of Theorem~2.3 also in this case, 
thus avoiding the assumption (A2).}
\hfill $\Box$
\medskip

Moreover, \pier{in the framework of assumptions \hbox{(A1)}, \hbox{(A3)} with \eqref{error2}, \hbox{(A5)} and \hbox{(A6)},  we can also improve the error estimate stated in Theorem~5.1.}

\paragraph{Theorem~6.1.} 
{\it Assume \hbox{(A1)}, \hbox{(A3)}, \hbox{(A5)} and \hbox{(A6)} with \eqref{error2} and \eqref{error3}. 
For $\varepsilon \in (0,\bar{\varepsilon }]$, let 
$(u_\varepsilon ,\mu_\varepsilon ,\xi _\varepsilon )$ be a solution of 
problem \hbox{(P)}$_\varepsilon $ and let $(u,\mu ,\xi )$ be a solution of 
the problem \hbox{(P)}. Then, there exists a constant $\pier{C^{\star}}>0$, 
depending only on the data, such that}
\begin{equation} 
	| u_\varepsilon -  u 
	|_{C([0,T];V^*)}^2 
	+ 
	\int_{0}^{T} \bigl( \xi _\varepsilon (s)-\xi (s), u_\varepsilon (s)-u(s) \bigr)_{\! H} ds 
	\le \pier{C^{\star}} \varepsilon ^{1/2}
	\label{error*}
\end{equation} 
{\it for all $\varepsilon \in (0,\bar{\varepsilon }]$. }

\paragraph{Proof.}
Compare  \pier{\eqref{boot4} and \eqref{6-6}: it turns out that  
the boundedness  in $L^2(0,T;H)$ of 
$\{ \varepsilon^{2/3} \Delta u_\varepsilon \}_{\varepsilon >0}$ 
is improved to  the one of 
$\{ \varepsilon^{1/2} \Delta u_\varepsilon \}_{\varepsilon >0}$.
Therefore, in the proof of Theorem~\pier{5.1} 
the key estimate \eqref{key1} can be  replaced by}
\begin{align} 
	& 
	\varepsilon \int_{0}^{T}
	\bigl( \Delta u_\varepsilon (s),u_\varepsilon (s)-u(s) \bigr)_{\! H} ds
	\nonumber \\
	& \quad  \le 
	\varepsilon ^{1/2} 
	\left\{ \varepsilon  \int_{0}^{T} \bigl| \Delta u_\varepsilon (s) \bigr|_H^2 ds \right\}^{1/2} 
	\left\{ \int_{0}^{T} \bigl| u_\varepsilon (s)-u(s) \bigr|_H^2 ds \right\}^{1/2} 
	 \le \pier{C_4^*} \varepsilon ^{1/2} ,
	\label{key*}
\end{align}
\pier{for some constant $C_4^*$  depending only on $c_1$, $c_3$, $c_4$, $ |\Delta f|_{L^2(0,T;H)}$, $M_1$, $|\Omega|$ and $T$. Recalling now \eqref{error3} and 
\eqref{key2}, it is straightforward to derive \eqref{error*} from \eqref{key}.}
\hfill $\Box$

\medskip
\pier{The convergence result and the error estimate shown in this section turn out to be an  
improvement which does without the assumption (A2) and exploits (A6) instead. Actually, 
thanks to this, we can treat more general cases of $\beta $, in particular 
they apply to the problems described in Examples~5 and~6.}

\paragraph{Remark~\pier{6.2}.} 
We note that in the framework of this section (cf.\ (A6)), the limit problem in 
\eqref{weak1}--\eqref{weak3} can be rewritten as
\begin{equation*} 
	\bigl \langle u'(t),z
	\bigr \rangle _{V^*,V}
	+ \int_{\Omega }^{} \nabla \xi(t) 
	\cdot \nabla z dx
	=  \int_{\Omega }^{} \pier{g(t)}\, z dx
	\quad \hbox{for~all~} z \in V,
\end{equation*}
for a.a.\ $t\in (0,T)$, with
\begin{gather*} 
	\xi  \in 
	\beta (u)
	\quad \hbox{a.e.~in~} Q, \quad 
	u(0)=\pier{u_0}
	\quad \hbox{a.e.~in~} \Omega.
\end{gather*} 
Here, the variable $\mu $ \pier{disappears. We point out that 
now it is also $\xi $, and not only $\xi -f$, to formally satisfy the {N}eumann 
homogeneous boundary condition 
\begin{equation} 
	\partial_{\boldsymbol{\nu }} \xi =0 
	\quad \hbox{a.e.\ on}~\Sigma.
\label{pier6}
\end{equation}
Owing to the zero mean value condition for $g$, now 
\eqref{pier6} is a necessary} condition for the existence of solutions.

\section{Appendix}
\setcounter{equation}{0}

We use the same notation as in the previous sections for function spaces.

\paragraph{Lemma A.} 
{\it \pier{Let} $u_0 \in H$ with 
$\widehat{\beta }(u_0) \in L^1(\Omega)$ and \pier{$m_0 \in \interior D(\beta )$}. Then there 
exist \pier{a family $\{u_{0\varepsilon }\}_{\varepsilon \in (0,1]}\subset  V$ and a positive constant $C$ such that}}
\begin{gather*} 
\pier{m(u_{0\varepsilon })=m_0,}\quad 
	|u_{0\varepsilon }|_H^2 \le C, 
	\quad 
	\int_{\Omega }^{} \widehat{\beta }(u_{0\varepsilon }) dx \le C, 
	\quad 
	\varepsilon |\nabla u_{0\varepsilon }|^2_{H^d} \le C
	\quad \hbox{\it for~all }\varepsilon \in (0,1],
	\\
	u_{0\varepsilon } \to u \quad \pier{\hbox{\it strongly~in~$ H $ as } \varepsilon \searrow 0.}
\end{gather*} 

\paragraph{Proof.} For each $\varepsilon \in (0,1]$, we can take $u_{0\varepsilon } \pier{{}\in W}$ as 
the solution of 
\begin{equation*} 
	\begin{cases}
	\displaystyle u_{0\varepsilon }-\varepsilon \Delta u_{0\varepsilon } = u_0 
	\quad \hbox{a.e.\ in~} \Omega, \\[2mm]
	\partial _{\boldsymbol{\nu }} u_{0\varepsilon } = 0 
	\quad \hbox{a.e.\ in~} \Gamma. \\
	\end{cases} 
\end{equation*}
Then, $m(u_{0\varepsilon })=m(u_0)=m_0$ and 
$u_{0\varepsilon } \to u_0$ strongly in $H$ as $\varepsilon \searrow 0$. 
Indeed, testing the first equation by $u_{0\varepsilon }$
and using the Young inequality, we find 
\begin{align}
	\int_{\Omega }^{} |u_{0\varepsilon }|^2 dx 
	+ \varepsilon \int_{\Omega }^{} 
	|\nabla u_{0\varepsilon }|^2 dx 
	& \le \pier{\frac{1}{2} \int_{\Omega }^{} |u_0|^2 dx 
	+ \frac{1}{2} \int_{\Omega }^{} |u_{0\varepsilon }|^2 dx,} 
	\label{7-1}
\end{align}
whence $\{ u_{0\varepsilon } \}_{\varepsilon \in (0,1]}$ is bounded in $H$ and 
$\varepsilon u_{0\varepsilon } \to 0$ strongly in $V$ as $\varepsilon \searrow 0$.  
Then, from the equation it turns out that 
$u_{0\varepsilon } \to u_0$ weakly in $H$, when passing to the limit 
in 
\begin{equation}
	\int_{\Omega }^{} 
	(u_{0\varepsilon }-u_0)z dx + \varepsilon 
	\int_{\Omega}^{} \nabla u_{0\varepsilon } \cdot \nabla z dx=0 
	\quad {\rm for~all}~z \in V.
	\label{apkey}
\end{equation}
Moreover, from \eqref{7-1} \pier{it follows that}
\begin{equation*}
	\limsup_{\varepsilon \searrow  0} 
	\int_{\Omega }^{} 
	|u_{0\varepsilon }|^2 dx \le 
	\int_{\Omega }^{}
	|u_0|^2 dx.
\end{equation*}
\pier{This ensures convergence of norms and finally strong convergence of
$u_{0\varepsilon } $ to $u_0$ in $H$.} 
Next, taking $z:=\beta_{\tilde{\varepsilon}} (u_{0\varepsilon })$, where 
$\beta _{\tilde{\varepsilon }}$ is the {Y}osida approximation of $\beta $ 
\pier{(treated in Section~3) at
$\tilde{\varepsilon } \in (0,1]$,} and 
using the definition of the subdifferential lead to 
\begin{align*}
	\int_{\Omega }^{} 
	\bigl(
	\widehat{\beta }_{\tilde{\varepsilon }}(u_{0\varepsilon })
	- \widehat{\beta }_{\tilde{\varepsilon }}(u_{0})
	\bigr)dx
	& \le \int_{\Omega }^{} (u_{0\varepsilon }-u_0) \beta _{\tilde{\varepsilon }}(u_{0\varepsilon })dx \\
	& \pier{{}= - \varepsilon  \int_{\Omega }^{} \beta '_{\tilde{\varepsilon }}(u_{0\varepsilon }) |\nabla u_{0\varepsilon }|^2 dx  \le 0.}
\end{align*}
\pier{Therefore, we have that} 
\begin{equation*} 
	\int_{\Omega }^{} 
	\widehat{\beta }_{\tilde{\varepsilon }}(u_{0\varepsilon }) dx 
	\le 
	\int_{\Omega }^{} 
	\widehat{\beta }_{\tilde{\varepsilon }}(u_0) dx 
	\le \int_{\Omega }^{} 
	\widehat{\beta }(u_0) dx \le C;
\end{equation*} 
if we take $\tilde{\varepsilon } <\varepsilon $ and pass to the limit \pier{as} $\tilde{\varepsilon } \searrow 0$,  then we deduce that 
\begin{equation*} 
	\int_{\Omega }^{} 
	\widehat{\beta }(u_{0\varepsilon }) dx 
	=  \lim _{\tilde{\varepsilon } \searrow  0}
	\int_{\Omega }^{} 
	\widehat{\beta }_{\tilde{\varepsilon }}(u_{0\varepsilon }) dx 
	\le C. 
\end{equation*} 
Thus, \eqref{apini} is completely proved. 
\pier{As a remark, the additional condition \eqref{error3} is also guaranteed since we can take}
$z:={\cal N}(u_{0\varepsilon }-u_0)$ in \eqref{apkey}. 
\hfill $\Box$

\section*{Acknowledgments}

The authors wish to express their heartfelt gratitude to
professors {G}oro {A}kagi and {U}lisse {S}tefanelli, 
who kindly gave them the opportunity of exchange visits, 
supported by the JSPS--CNR bilateral joint research 
\emph{{I}nnovative {V}ariational {M}ethods for {E}volution 
{E}quations}. \pier{The present paper also benefits from the support of the MIUR-PRIN Grant 2010A2TFX2 ``Calculus of Variations'' and the GNAMPA (Gruppo Nazionale per l'Analisi Matematica, la Probabilit\`a e le loro Applicazioni) of INdAM (Istituto Nazionale di Alta Matematica) for PC, and of the JSPS KAKENHI 
Grant-in-Aid for Scientific Research(C), Grant Number 26400164 for TF.}

\end{document}